\documentclass[reqno,12pt,a4paper]{amsart}

\usepackage{mathrsfs}
\usepackage{amssymb}
\usepackage{amsfonts}
\usepackage{latexsym}
\usepackage{amsthm}
\usepackage{graphicx}
\usepackage{xcolor}

\def\lb{\label}

\newcommand{\er}[1]{\textrm{(\ref{#1})}}

\begin{document}


\renewcommand{\theequation}{\arabic{section}.\arabic{equation}}
\theoremstyle{plain}
\newtheorem{theorem}{\bf Theorem}[section]
\newtheorem{lemma}[theorem]{\bf Lemma}
\newtheorem{corollary}[theorem]{\bf Corollary}
\newtheorem{proposition}[theorem]{\bf Proposition}
\newtheorem{definition}[theorem]{\bf Definition}

\theoremstyle{remark}
\newtheorem{remark}[theorem]{\bf Remark}
\newtheorem{example}[theorem]{\bf Example}

\def\a{\alpha}  \def\cA{{\mathcal A}}     \def\bA{{\bf A}}  \def\mA{{\mathscr A}}
\def\b{\beta}   \def\cB{{\mathcal B}}     \def\bB{{\bf B}}  \def\mB{{\mathscr B}}
\def\g{\gamma}  \def\cC{{\mathcal C}}     \def\bC{{\bf C}}  \def\mC{{\mathscr C}}
\def\G{\Gamma}  \def\cD{{\mathcal D}}     \def\bD{{\bf D}}  \def\mD{{\mathscr D}}
\def\d{\delta}  \def\cE{{\mathcal E}}     \def\bE{{\bf E}}  \def\mE{{\mathscr E}}
\def\D{\Delta}  \def\cF{{\mathcal F}}     \def\bF{{\bf F}}  \def\mF{{\mathscr F}}
\def\c{\chi}    \def\cG{{\mathcal G}}     \def\bG{{\bf G}}  \def\mG{{\mathscr G}}
\def\z{\zeta}   \def\cH{{\mathcal H}}     \def\bH{{\bf H}}  \def\mH{{\mathscr H}}
\def\e{\eta}    \def\cI{{\mathcal I}}     \def\bI{{\bf I}}  \def\mI{{\mathscr I}}
\def\p{\psi}    \def\cJ{{\mathcal J}}     \def\bJ{{\bf J}}  \def\mJ{{\mathscr J}}
\def\vT{\Theta} \def\cK{{\mathcal K}}     \def\bK{{\bf K}}  \def\mK{{\mathscr K}}
\def\k{\kappa}  \def\cL{{\mathcal L}}     \def\bL{{\bf L}}  \def\mL{{\mathscr L}}
\def\l{\lambda} \def\cM{{\mathcal M}}     \def\bM{{\bf M}}  \def\mM{{\mathscr M}}
\def\L{\Lambda} \def\cN{{\mathcal N}}     \def\bN{{\bf N}}  \def\mN{{\mathscr N}}
\def\m{\mu}     \def\cO{{\mathcal O}}     \def\bO{{\bf O}}  \def\mO{{\mathscr O}}
\def\n{\nu}     \def\cP{{\mathcal P}}     \def\bP{{\bf P}}  \def\mP{{\mathscr P}}
\def\r{\rho}    \def\cQ{{\mathcal Q}}     \def\bQ{{\bf Q}}  \def\mQ{{\mathscr Q}}
\def\s{\sigma}  \def\cR{{\mathcal R}}     \def\bR{{\bf R}}  \def\mR{{\mathscr R}}
\def\S{\Sigma}  \def\cS{{\mathcal S}}     \def\bS{{\bf S}}  \def\mS{{\mathscr S}}
\def\t{\tau}    \def\cT{{\mathcal T}}     \def\bT{{\bf T}}  \def\mT{{\mathscr T}}
\def\f{\phi}    \def\cU{{\mathcal U}}     \def\bU{{\bf U}}  \def\mU{{\mathscr U}}
\def\F{\Phi}    \def\cV{{\mathcal V}}     \def\bV{{\bf V}}  \def\mV{{\mathscr V}}
\def\P{\Psi}    \def\cW{{\mathcal W}}     \def\bW{{\bf W}}  \def\mW{{\mathscr W}}
\def\o{\omega}  \def\cX{{\mathcal X}}     \def\bX{{\bf X}}  \def\mX{{\mathscr X}}
\def\x{\xi}     \def\cY{{\mathcal Y}}     \def\bY{{\bf Y}}  \def\mY{{\mathscr Y}}
\def\X{\Xi}     \def\cZ{{\mathcal Z}}     \def\bZ{{\bf Z}}  \def\mZ{{\mathscr Z}}
\def\be{{\bf e}} \def\bc{{\bf c}}
\def\bv{{\bf v}} \def\bu{{\bf u}}
\def\Om{\Omega}
\def\bp{{\bf p}}\def\bq{{\bf q}}
\def\bx{{\bf x}} \def\by{{\bf y}}
\def\bbD{\pmb \Delta}
\def\mm{\mathrm m}
\def\mn{\mathrm n}

\newcommand{\mc}{\mathscr {c}}

\newcommand{\gA}{\mathfrak{A}}          \newcommand{\ga}{\mathfrak{a}}
\newcommand{\gB}{\mathfrak{B}}          \newcommand{\gb}{\mathfrak{b}}
\newcommand{\gC}{\mathfrak{C}}          \newcommand{\gc}{\mathfrak{c}}
\newcommand{\gD}{\mathfrak{D}}          \newcommand{\gd}{\mathfrak{d}}
\newcommand{\gE}{\mathfrak{E}}
\newcommand{\gF}{\mathfrak{F}}           \newcommand{\gf}{\mathfrak{f}}
\newcommand{\gG}{\mathfrak{G}}           
\newcommand{\gH}{\mathfrak{H}}           \newcommand{\gh}{\mathfrak{h}}
\newcommand{\gI}{\mathfrak{I}}           \newcommand{\gi}{\mathfrak{i}}
\newcommand{\gJ}{\mathfrak{J}}           \newcommand{\gj}{\mathfrak{j}}
\newcommand{\gK}{\mathfrak{K}}            \newcommand{\gk}{\mathfrak{k}}
\newcommand{\gL}{\mathfrak{L}}            \newcommand{\gl}{\mathfrak{l}}
\newcommand{\gM}{\mathfrak{M}}            \newcommand{\gm}{\mathfrak{m}}
\newcommand{\gN}{\mathfrak{N}}            \newcommand{\gn}{\mathfrak{n}}
\newcommand{\gO}{\mathfrak{O}}
\newcommand{\gP}{\mathfrak{P}}             \newcommand{\gp}{\mathfrak{p}}
\newcommand{\gQ}{\mathfrak{Q}}             \newcommand{\gq}{\mathfrak{q}}
\newcommand{\gR}{\mathfrak{R}}             \newcommand{\gr}{\mathfrak{r}}
\newcommand{\gS}{\mathfrak{S}}              \newcommand{\gs}{\mathfrak{s}}
\newcommand{\gT}{\mathfrak{T}}             \newcommand{\gt}{\mathfrak{t}}
\newcommand{\gU}{\mathfrak{U}}             \newcommand{\gu}{\mathfrak{u}}
\newcommand{\gV}{\mathfrak{V}}             \newcommand{\gv}{\mathfrak{v}}
\newcommand{\gW}{\mathfrak{W}}             \newcommand{\gw}{\mathfrak{w}}
\newcommand{\gX}{\mathfrak{X}}               \newcommand{\gx}{\mathfrak{x}}
\newcommand{\gY}{\mathfrak{Y}}              \newcommand{\gy}{\mathfrak{y}}
\newcommand{\gZ}{\mathfrak{Z}}             \newcommand{\gz}{\mathfrak{z}}

\def\ve{\varepsilon}   \def\vt{\vartheta}    \def\vp{\varphi}    \def\vk{\varkappa} \def\vr{\varrho}

\def\A{{\mathbb A}} \def\B{{\mathbb B}} \def\C{{\mathbb C}}
\def\dD{{\mathbb D}} \def\E{{\mathbb E}} \def\dF{{\mathbb F}} \def\dG{{\mathbb G}} \def\H{{\mathbb H}}\def\I{{\mathbb I}} \def\J{{\mathbb J}} \def\K{{\mathbb K}} \def\dL{{\mathbb L}}\def\M{{\mathbb M}} \def\N{{\mathbb N}} \def\O{{\mathbb O}} \def\dP{{\mathbb P}} \def\R{{\mathbb R}}\def\S{{\mathbb S}} \def\T{{\mathbb T}} \def\U{{\mathbb U}} \def\V{{\mathbb V}}\def\W{{\mathbb W}} \def\X{{\mathbb X}} \def\Y{{\mathbb Y}} \def\Z{{\mathbb Z}}


\def\la{\leftarrow}              \def\ra{\rightarrow}            \def\Ra{\Rightarrow}
\def\ua{\uparrow}                \def\da{\downarrow}
\def\lra{\leftrightarrow}        \def\Lra{\Leftrightarrow}


\def\lt{\biggl}                  \def\rt{\biggr}
\def\ol{\overline}               \def\wt{\widetilde}
\def\ul{\underline}
\def\no{\noindent}


\let\ge\geqslant                 \let\le\leqslant
\def\lan{\langle}                \def\ran{\rangle}
\def\/{\over}                    \def\iy{\infty}
\def\sm{\setminus}               \def\es{\emptyset}
\def\ss{\subset}                 \def\ts{\times}
\def\pa{\partial}                \def\os{\oplus}
\def\om{\ominus}                 \def\ev{\equiv}
\def\iint{\int\!\!\!\int}        \def\iintt{\mathop{\int\!\!\int\!\!\dots\!\!\int}\limits}
\def\el2{\ell^{\,2}}             \def\1{1\!\!1}
\def\sh{\sharp}
\def\wh{\widehat}
\def\bs{\backslash}
\def\intl{\int\limits}

\def\na{\mathop{\mathrm{\nabla}}\nolimits}
\def\sh{\mathop{\mathrm{sh}}\nolimits}
\def\ch{\mathop{\mathrm{ch}}\nolimits}
\def\where{\mathop{\mathrm{where}}\nolimits}
\def\all{\mathop{\mathrm{all}}\nolimits}
\def\as{\mathop{\mathrm{as}}\nolimits}
\def\Area{\mathop{\mathrm{Area}}\nolimits}
\def\arg{\mathop{\mathrm{arg}}\nolimits}
\def\const{\mathop{\mathrm{const}}\nolimits}
\def\det{\mathop{\mathrm{det}}\nolimits}
\def\diag{\mathop{\mathrm{diag}}\nolimits}
\def\diam{\mathop{\mathrm{diam}}\nolimits}
\def\dim{\mathop{\mathrm{dim}}\nolimits}
\def\dist{\mathop{\mathrm{dist}}\nolimits}
\def\Im{\mathop{\mathrm{Im}}\nolimits}
\def\Iso{\mathop{\mathrm{Iso}}\nolimits}
\def\Ker{\mathop{\mathrm{Ker}}\nolimits}
\def\Lip{\mathop{\mathrm{Lip}}\nolimits}
\def\rank{\mathop{\mathrm{rank}}\limits}
\def\Ran{\mathop{\mathrm{Ran}}\nolimits}
\def\Re{\mathop{\mathrm{Re}}\nolimits}
\def\Res{\mathop{\mathrm{Res}}\nolimits}
\def\res{\mathop{\mathrm{res}}\limits}
\def\sign{\mathop{\mathrm{sign}}\nolimits}
\def\span{\mathop{\mathrm{span}}\nolimits}
\def\supp{\mathop{\mathrm{supp}}\nolimits}
\def\Tr{\mathop{\mathrm{Tr}}\nolimits}
\def\BBox{\hspace{1mm}\vrule height6pt width5.5pt depth0pt \hspace{6pt}}


\newcommand\nh[2]{\widehat{#1}\vphantom{#1}^{(#2)}}
\def\dia{\diamond}

\def\Oplus{\bigoplus\nolimits}



\def\qqq{\qquad}
\def\qq{\quad}
\let\ge\geqslant
\let\le\leqslant
\let\geq\geqslant
\let\leq\leqslant
\newcommand{\ca}{\begin{cases}}
\newcommand{\ac}{\end{cases}}
\newcommand{\ma}{\begin{pmatrix}}
\newcommand{\am}{\end{pmatrix}}
\renewcommand{\[}{\begin{equation}}
\renewcommand{\]}{\end{equation}}
\def\eq{\begin{equation}}
\def\qe{\end{equation}}
\def\[{\begin{equation}}
\def\bu{\bullet}

\makeatletter
\@namedef{subjclassname@2020}{\textup{2020} Mathematics Subject Classification}
\makeatother

\title{Spectrum of Schr\"odinger operators on subcovering graphs}

\date{\today}
\author[Natalia Saburova]{Natalia Saburova}
\address{Northern (Arctic) Federal University, Severnaya Dvina Emb. 17, Arkhangelsk, 163002, Russia, \ n.saburova@gmail.com, \ n.saburova@narfu.ru}


\begin{abstract}
We consider discrete Schr\"odinger operators with real periodic potentials on periodic graphs. The spectra of the operators consist of a finite number of bands. By "rolling up" a periodic graph along some appropriate directions we obtain periodic graphs of smaller dimensions called subcovering graphs. For example, rolling up a planar hexagonal lattice along different directions will lead to nanotubes with various chiralities. We describe connections between spectra of the Schr\"odinger operators on a periodic graph and its subcoverings. In particular, we provide a simple criterion for the subcovering graph to be isospectral to the original periodic graph. By isospectrality of periodic graphs we mean that the spectra of the Schr\"odinger operators on the graphs consist of the same number of bands and the corresponding bands coincide as sets. We also obtain asymptotics of the band edges of the Schr\"odinger operator on the subcovering graph as the "chiral" (roll up) vectors are long enough.
\end{abstract}

\maketitle

\section {\lb{Sec0}{Introduction}}
\setcounter{equation}{0}

Schr\"odinger operators on graphs have a lot of applications in physics, chemistry and engineering.  In solid-state physics, the tight-binding model is a commonly used approach to the calculation of electronic band structure of solids \cite{AM76}. The material is modeled by a discrete graph with vertices at the location of atoms in solid and with edges indicating chemical bonding of atoms.

In this paper we consider discrete Schr\"odinger operators with real periodic potentials on arbitrary periodic graphs. Their spectra consist of a finite number of segments called \emph{bands}. The $n$-th band is the range of a \emph{band} function $\l_n(k)$ describing the dependence of electron energy on the multidimensional parameter $k$ called \emph{quasimomentum} which takes its values in the \emph{Brillouin zone} $(-\pi,\pi]^d$, where $d$ is the dimension of the periodic graph (i.e., the rank of its period lattice). The multi-valued function $k\mapsto\{\l_n(k)\}$ is called the \emph{dispersion relation}. By "rolling up" a periodic graph in some appropriate directions we obtain periodic graphs of smaller dimensions called subcovering graphs. The aim of this paper is to describe connections between spectra of the Schr\"odinger operators on a periodic graph and its subcoverings.

This paper was partly motivated by the nice article \cite{KP07} describing simple connections between spectra of Schr\"odinger operators on graphene and nanotubes. Graphene is a single 2D layer of graphite forming a hexagonal lattice. The hexagonal lattice is invariant under translations by vectors $\ga_1$ and $\ga_2$, see Fig.~\ref{fig1}\emph{a}. A carbon nanotube is a hexagonal lattice "rolled up" into a cylinder. In carbon nanotubes, the graphene sheet is "rolled up" in such a way that a so-called \emph{chiral} vector $\gt=t_1\ga_1+t_2\ga_2$, where $t:=(t_1,t_2)\in\Z^2$, becomes the circumference of the tube, see Fig.~\ref{fig1}\emph{b,c}. In \cite{KP07}, applying a simple restriction procedure to the dispersion relation for the Schr\"odinger operator on the graphene lattice, the authors described the spectra of the Schr\"odinger operators on any carbon nanotube: zigzag ($t=(n,0)$), armchair ($t=(n,n)$) or chiral. In particular, they showed that the spectra of the Schr\"odinger operators on a nanotube and graphene coincide as sets if and only if $t_1-t_2$ is divisible by 3, and these spectra are asymptotically close as the length of the chiral vector $\gt$ tends to infinity. The similar approach was used in \cite{D15} to derive spectral properties of \emph{graphyne} nanotubes from the dispersion relation of the \emph{graphyne} lattice. In \cite{KS21} the authors also treated the connection between spectra of the Laplacians on a $\Z^2$-periodic 3-regular graph and its subcoverings to construct graphs with maximal gap intervals.

Each periodic graph is a subcovering of some higher dimensional periodic graphs and, in particular, of a so-called maximal abelian covering graph (see Chapter 6 in \cite{S13}). Sometimes, due to the presence of some symmetries, it is easier to describe the spectrum of the operators on such higher dimensional periodic graphs than on their subcoverings. For example, there are some known results about spectral properties of operators on the maximal abelian covering graphs. In \cite{HN09,HS99,HS04} the maximal abelian coverings on which the spectra of Laplacians have no gaps and the maximal abelian coverings with no eigenvalues (degenerate bands) were described. A more general class of periodic graphs with the minimal number of links between copies of the fundamental domain ("one crossing edge per generator") was considered in \cite{BCCM22}. In particular, this class of periodic graphs includes all maximal abelian coverings and their "pendant" and "spider" decorations (see, e.g., \cite{K05}). In \cite{BCCM22} it was shown that any local extremum of a band function of the discrete Schr\"odinger operators on such graphs is in fact its global extremum if the dimension $d$ of the graph is at most three, or (in any dimension) if the critical point is a \emph{symmetry} point of the Brillouin zone from the set $\{0,\pi\}^d$. For the discussion of the question why the band edges for periodic operators are often (but not always!) attained at the symmetry points of the Brillouin zone see also \cite{HKSW07}.

Studying connections between spectra of operators on a periodic graph and its subcovering graphs allows one to derive some spectral properties of subcoverings from the spectrum of the operator on the original higher dimensional periodic graph (for example, on the corresponding maximal abelian covering graph or periodic graph with "one crossing edge per generator") which might be determined more easily.

In this paper we describe spectra of the Schr\"odinger operators on subcovering graphs using some apriori information about the spectrum of the operator on the original periodic graph. In particular, we formulate a criterion for the subcovering graph to be isospectral to the original periodic graph. By isospectrality of periodic graphs we mean that the spectra of the Schr\"odinger operators on the graphs consist of the same number of bands and the corresponding bands coincide as sets. We get asymptotics of the band edges of the Schr\"odinger operator on the subcoverings as the length of the chiral (roll up) vectors tends to infinity. The proof is based on the restriction procedure applied to the dispersion relation for the Schr\"odinger operator on the original periodic graph which was used in \cite{KP07} for the hexagonal lattice, and on the construction of the minimum norm solution of a linear system of equations corresponding to this restriction.

We also mention the recent paper \cite{S24}, where a specific case of subcovering graphs was studied and the notion of asymptotic isospectrality for periodic graphs was introduced.  Namely, in \cite{S24} the author considered perturbations of a periodic graph by adding edges in a periodic way (without changing the vertex set) and showed that if the added edges are long enough, then the perturbed graph is asymptotically isospectral to some periodic graph of a higher dimension but without long edges (called a \emph{limit} graph). This limit graph is actually an infinite-fold covering graph over the perturbed one.

\subsection{Discrete Schr\"odinger operators on periodic graphs.}
Let $\cG=(\cV,\cE)$ be a connected infinite graph embedded into $\R^d$. Here
$\cV$  is the set of its vertices and $\cE$ is the set of its (unoriented) edges. An edge between vertices $u,v\in\cV$ will be denoted as the (unordered) pair $\{u,v\}$ and said to be \emph{incident} to the vertices $u$ and $v$.

Let $\G$ be a lattice of rank $d$ in $\R^d$ with a basis $\{\ga_1,\ldots,\ga_d\}$, i.e.,
$$
\G=\Big\{\ga\in\R^d : \ga=\sum_{s=1}^dn_s\ga_s, \; (n_s)_{s=1}^d\in\Z^d\Big\},
$$
and let
$$
\Omega=\Big\{\bx\in\R^d : \bx=\sum_{s=1}^dx_s\ga_s, \; (x_s)_{s=1}^d\in
[0,1)^d\Big\}
$$
be the \emph{fundamental cell} of the lattice $\G$. We define the equivalence relation on $\R^d$:
$$
\bx\equiv \by \; (\hspace{-4mm}\mod \G) \qq\Leftrightarrow\qq \bx-\by\in\G \qqq
\forall\, \bx,\by\in\R^d.
$$

We consider \emph{locally finite $\G$-periodic graphs} $\cG$, i.e.,
graphs satisfying the following conditions:
\begin{itemize}
  \item[1)] $\cG=\cG+\ga$ for any $\ga\in\G$, i.e., $\cG$ is invariant under translation by any vector $\ga\in\G$;
  \item[2)] the quotient graph  $\cG_*=\cG/\G$ is finite.
\end{itemize}
The basis vectors $\ga_1,\ldots,\ga_d$ of the lattice $\G$ are called the {\it periods}  of $\cG$. We also call the quotient graph $\cG_*=\cG/\G$ the \emph{fundamental graph} of the periodic graph $\cG$. The fundamental graph $\cG_*=(\cV_*,\cE_*)$ has the vertex set $\cV_*=\cV/\G$ and the edge set $\cE_*=\cE/\G$ which are finite.

Let $\ell^2(\cV)$ be the Hilbert space of all square summable
functions  $f:\cV\to \C$ equipped with the norm
$$
\|f\|^2_{\ell^2(\cV)}=\sum_{v\in\cV}|f(v)|^2<\iy.
$$
On edges of $\cG=(\cV,\cE)$ we define a positive $\G$-periodic weight
$$
\o:\cE\to(0,+\iy),\qqq \o(\be+\ga)=\o(\be),\qqq \forall\,(\be,\ga)\in\cE\ts\G,
$$
and consider the Schr\"odinger operator $H$ acting on $\ell^2(\cV)$ and given by
\[\lb{Sh}
H=\D+Q.
\]
Here $\D$ is the discrete weighted Laplacian having the form
\[\lb{ALO}
(\D f)(v)=\sum_{\be=\{v,u\}\in\cE}\o(\be)\big(f(v)-f(u)\big), \qqq f\in\ell^2(\cV), \qqq v\in\cV,
\]
and $Q$ is a real $\G$-periodic potential, i.e., it satisfies
$$
(Qf)(v)=Q(v)f(v), \qqq Q(v+\ga)=Q(v), \qqq \forall\,(v,\ga)\in\cV\ts\G.
$$
The sum in (\ref{ALO}) is taken over all edges $\be\in\cE$ incident to the vertex $v$. It is known that the Schr\"odinger operator $H$ is a bounded self-adjoint operator on $\ell^2(\cV)$ (see, e.g., \cite{SS92}).

We note that the embedding of the periodic graph $\cG$ into $\R^d$ plays no role in the analysis of $H$, since only the incidence structure of $\cG$ matters. Thus, without loss of generality we may assume that \textbf{the period lattice $\G$ of $\cG$ is just the integer lattice $\Z^d$}.

To describe the spectrum of the Schr\"odinger operator $H$ on periodic graphs one can use the standard Floquet-Bloch theory (see, e.g., \cite{RS78} or for the graph case \cite{BK13}, Chapter 4). We introduce a family of spaces $\ell^2_k$ depending on the parameter $k\in\cB$ called \emph{quasimomentum}, where $\cB=(-\pi,\pi]^d$ is the \emph{Brillouin zone}. For each fixed $k$ the space $\ell^2_k$ consists of all functions $f:\cV\to \C$ satisfying the \emph{Floquet-Bloch} condition
\[\lb{FlBc}
f(v+\ga)=e^{i\lan \ga,k\ran}f(v),\qqq \forall\, v\in\cV,\qqq\forall\,\ga\in\Z^d,
\]
where $\lan \cdot,\cdot\ran$ denotes the standard inner product in $\R^d$. Such a function $f$ is uniquely determined by its values at the vertices of the periodic graph from the fundamental cell $\Omega$, and thus $\ell^2_k$ is naturally isomorphic to $\ell^2(\cV_*)$, where $\cV_*$ is the set of the fundamental graph vertices. We define the \emph{Floquet operator} as follows
\[\lb{Flop}
H(k):\ell^2(\cV_*)\to\ell^2(\cV_*)\qq\textrm{is the restriction of $H$ to $\ell^2_k$.}
\]

\begin{remark}\lb{Rfo0}
The Floquet operator $H(0)$ is just the Schr\"odinger operator defined by \er{Sh}, \er{ALO} on the fundamental graph $\cG_*=(\cV_*,\cE_*)$.
\end{remark}

The Floquet-Bloch theory provides the direct integral expansion
\[\lb{die0}
H=\int^\oplus_{\cB}H(k)\,dk,\qqq \cB=(-\pi,\pi]^d,
\]
which, in particular, yields
$$
\s(H)=\bigcup_{k\in\cB}\s\big(H(k)\big).
$$

Each Floquet operator $H(k)$ is self-adjoint and has $\n=\#\cV_*$ real eigenvalues
$$
\l_1(k)\leq\l_2(k)\leq\ldots\leq\l_\n(k), \qqq k\in\cB,
$$
labeled in non-decreasing order counting multiplicities. Here $\# M$ denotes the number of elements in a set $M$. The multi-valued function $k\mapsto\{\l_j(k)\}_{j=1}^\n$ is called the \emph{dispersion relation} for the Schr\"odinger operator $H$. Each \emph{band function} $\l_j(\cdot)$ is a continuous and piecewise real analytic function on the torus $\R^d/2\pi\Z^d$ and creates the \emph{spectral band} $\s_j(H)$ given by
\[\lb{ban.1H}
\begin{array}{l}
\s_j(H)=[\l_j^-,\l_j^+]=\l_j(\cB),\qqq j\in\N_\n, \qqq \N_\n=\{1,\ldots,\n\},\\[8pt]
\displaystyle\textrm{where}\qqq \l_j^-=\min_{k\in\cB}\l_j(k),\qqq \l_j^+=\max_{k\in\cB}\l_j(k).
\end{array}
\]
Some of $\l_j(\cdot)$ may be constant, i.e., $\l_j(\cdot)=\L=\const$, on some subset of $\cB$ of positive Lebesgue measure. In this case the Schr\"odinger operator $H$ on $\cG$ has the eigenvalue $\L$ of infinite multiplicity. We
call $\L$ a \emph{flat band}. Thus, the spectrum of the Schr\"odinger operator $H$ on a periodic graph $\cG$ has the form (for more details see, e.g., \cite{HN09})
$$
\s(H)=\bigcup_{k\in\cB}\s\big(H(k)\big)=
\bigcup_{j=1}^\n\s_j(H)=\s_{ac}(H)\cup \s_{fb}(H),
$$
where $\s_{ac}(H)$ is the absolutely continuous spectrum, which is a
union of non-degenerate bands from \er{ban.1H}, and $\s_{fb}(H)$ is
the set of all flat bands.

\begin{remark}\lb{prbf}
\emph{i}) The first spectral band $\s_1(H)=[\l_1^-,\l_1^+]$ is never a flat band and $\l_1^{-}=\l_1(0)$, see Theorem 2.7 in \cite{SY23} and Theorem 1 in \cite{SS92}.

\emph{ii}) The band functions $\l_j(\cdot)$, $j\in\N_\n$, are even (see Lemma 2.9.\emph{i} in \cite{KKR17}), i.e.,
$$
\l_j(-k)=\l_j(k),\qqq \forall\,k\in\cB.
$$
\end{remark}

\subsection{Subcovering graphs}\lb{Sscg}
Let $\cG$ be a $\Z^d$-periodic graph. (Recall that we identify the period lattice $\G$ of $\cG$ with the integer lattice $\Z^d$). We fix some basis $\ga_1,\ldots,\ga_d$ of $\Z^d$ (the periods of the graph $\cG$). Below \textbf{the coordinates of all vectors of $\Z^d$ will be given with respect to this fixed basis $\ga_1,\ldots,\ga_d$}.

Let $\G_\gt$ be a sublattice of $\Z^d$ with a basis $\gt:=\{\gt_1,\ldots,\gt_{d_o}\}\ss\Z^d$, where $d_o<d$. We consider a \emph{subcovering} graph $\cG_\gt$ of a $\Z^d$-periodic graph $\cG$ obtained as the quotient of $\cG$ with respect to the sublattice $\G_\gt$, i.e.,
$$
\cG_\gt=\cG/\G_\gt.
$$
The subcovering graph $\cG_\gt$ is obtained from the periodic graph $\cG$ by the identification of points $\bx,\by\in\cG$ such that $\bx-\by$ belongs to the lattice $\G_\gt$. The subcovering $\cG_\gt$ is naturally embedded into the cylinder $\R^d/\G_\gt$. Following the terminology for nanotubes, we will call the vectors $\gt_1,\ldots,\gt_{d_o}$ the \emph{chiral} vectors of the subcovering $\cG_\gt$.

\begin{figure}[t!]\centering
\unitlength 1.0mm 
\linethickness{0.4pt}

\begin{picture}(50,53)(0,0)
\put(0,0){\emph{a})}
\put(0,51){$\bG$}

\put(9,3){\circle*{1.5}}
\put(21,3){\circle*{1.5}}
\put(33,3){\circle*{1.5}}
\put(45,3){\circle*{1.5}}

\put(3,6){\circle*{1.5}}
\put(15,6){\circle*{1.5}}
\put(27,6){\circle*{1.5}}
\put(39,6){\circle*{1.5}}
\put(51,6){\circle*{1.5}}

\put(3,12){\circle*{1.5}}
\put(15,12){\circle*{1.5}}
\put(27,12){\circle*{1.5}}
\put(39,12){\circle*{1.5}}
\put(51,12){\circle*{1.5}}

\put(9,15){\circle*{1.5}}
\put(21,15){\circle*{1.5}}
\put(33,15){\circle*{1.5}}
\put(45,15){\circle*{1.5}}
\put(9,21){\circle*{1.5}}
\put(21,21){\circle*{1.5}}
\put(33,21){\circle*{1.5}}
\put(45,21){\circle*{1.5}}

\put(3,24){\circle*{1.5}}
\put(15,24){\circle*{1.5}}
\put(27,24){\circle*{1.5}}
\put(39,24){\circle*{1.5}}
\put(51,24){\circle*{1.5}}

\put(3,30){\circle*{1.5}}
\put(15,30){\circle*{1.5}}
\put(27,30){\circle*{1.5}}
\put(39,30){\circle*{1.5}}
\put(51,30){\circle*{1.5}}

\put(9,33){\circle*{1.5}}
\put(21,33){\circle*{1.5}}
\put(33,33){\circle*{1.5}}
\put(45,33){\circle*{1.5}}
\put(9,39){\circle*{1.5}}
\put(21,39){\circle*{1.5}}
\put(33,39){\circle*{1.5}}
\put(45,39){\circle*{1.5}}

\put(3,42){\circle*{1.5}}
\put(15,42){\circle*{1.5}}
\put(27,42){\circle*{1.5}}
\put(39,42){\circle*{1.5}}
\put(51,42){\circle*{1.5}}

\put(3,48){\circle*{1.5}}
\put(15,48){\circle*{1.5}}
\put(27,48){\circle*{1.5}}
\put(39,48){\circle*{1.5}}
\put(51,48){\circle*{1.5}}

\put(9,51){\circle*{1.5}}
\put(21,51){\circle*{1.5}}
\put(33,51){\circle*{1.5}}
\put(45,51){\circle*{1.5}}
\put(9,3){\line(0,-1){3.00}}
\put(21,3){\line(0,-1){3.00}}
\put(33,3){\line(0,-1){3.00}}
\put(45,3){\line(0,-1){3.00}}

\put(9,51){\line(0,1){3.00}}
\put(21,51){\line(0,1){3.00}}
\put(33,51){\line(0,1){3.00}}
\put(45,51){\line(0,1){3.00}}

\bezier{60}(3,6)(1.5,5.25)(0,4.5)
\bezier{60}(3,12)(1.5,12.75)(0,13.5)
\bezier{60}(3,24)(1.5,23.25)(0,22.5)
\bezier{60}(3,30)(1.5,30.75)(0,31.5)
\bezier{60}(3,42)(1.5,41.25)(0,40.5)
\bezier{60}(3,48)(1.5,48.75)(0,49.5)

\bezier{60}(51,6)(52.5,5.25)(54,4.5)
\bezier{60}(51,12)(52.5,12.75)(54,13.5)
\bezier{60}(51,24)(52.5,23.25)(54,22.5)
\bezier{60}(51,30)(52.5,30.75)(54,31.5)
\bezier{60}(51,42)(52.5,41.25)(54,40.5)
\bezier{60}(51,48)(52.5,48.75)(54,49.5)

\put(9,3){\line(-2,1){6.00}}
\put(9,3){\line(2,1){6.00}}
\put(21,3){\line(-2,1){6.00}}
\put(21,3){\line(2,1){6.00}}
\put(33,3){\line(-2,1){6.00}}
\put(33,3){\line(2,1){6.00}}
\put(45,3){\line(-2,1){6.00}}
\put(45,3){\line(2,1){6.00}}
\put(3,6){\line(0,1){6.00}}
\put(15,6){\line(0,1){6.00}}
\put(27,6){\line(0,1){6.00}}
\put(39,6){\line(0,1){6.00}}
\put(51,6){\line(0,1){6.00}}
\put(9,15){\line(-2,-1){6.00}}
\put(9,15){\line(2,-1){6.00}}
\put(21,15){\line(-2,-1){6.00}}
\put(21,15){\line(2,-1){6.00}}
\put(33,15){\line(-2,-1){6.00}}
\put(33,15){\line(2,-1){6.00}}
\put(45,15){\line(-2,-1){6.00}}
\put(45,15){\line(2,-1){6.00}}
\put(9,15){\line(0,1){6.00}}
\put(21,15){\line(0,1){6.00}}
\put(33,15){\line(0,1){6.00}}
\put(45,15){\line(0,1){6.00}}
\put(9,21){\line(-2,1){6.00}}
\put(9,21){\line(2,1){6.00}}
\put(21,21){\line(-2,1){6.00}}
\put(21,21){\line(2,1){6.00}}
\put(33,21){\line(-2,1){6.00}}
\put(33,21){\line(2,1){6.00}}
\put(45,21){\line(-2,1){6.00}}
\put(45,21){\line(2,1){6.00}}
\put(3,24){\line(0,1){6.00}}
\put(15,24){\line(0,1){6.00}}
\put(27,24){\line(0,1){6.00}}
\put(39,24){\line(0,1){6.00}}
\put(51,24){\line(0,1){6.00}}
\put(9,33){\line(-2,-1){6.00}}
\put(9,33){\line(2,-1){6.00}}
\put(21,33){\line(-2,-1){6.00}}
\put(21,33){\line(2,-1){6.00}}
\put(33,33){\line(-2,-1){6.00}}
\put(33,33){\line(2,-1){6.00}}
\put(45,33){\line(-2,-1){6.00}}
\put(45,33){\line(2,-1){6.00}}
\put(9,33){\line(0,1){6.00}}
\put(21,33){\line(0,1){6.00}}
\put(33,33){\line(0,1){6.00}}
\put(45,33){\line(0,1){6.00}}
\put(9,39){\line(-2,1){6.00}}
\put(9,39){\line(2,1){6.00}}
\put(21,39){\line(-2,1){6.00}}
\put(21,39){\line(2,1){6.00}}
\put(33,39){\line(-2,1){6.00}}
\put(33,39){\line(2,1){6.00}}
\put(45,39){\line(-2,1){6.00}}
\put(45,39){\line(2,1){6.00}}
\put(3,42){\line(0,1){6.00}}
\put(15,42){\line(0,1){6.00}}
\put(27,42){\line(0,1){6.00}}
\put(39,42){\line(0,1){6.00}}
\put(51,42){\line(0,1){6.00}}
\put(9,51){\line(-2,-1){6.00}}
\put(9,51){\line(2,-1){6.00}}
\put(21,51){\line(-2,-1){6.00}}
\put(21,51){\line(2,-1){6.00}}
\put(33,51){\line(-2,-1){6.00}}
\put(33,51){\line(2,-1){6.00}}
\put(45,51){\line(-2,-1){6.00}}
\put(45,51){\line(2,-1){6.00}}

\put(21.8,15.5){\small$v_1$}
\put(10.5,10.3){\small$v_2$}

\put(26,17){$\Omega$}

\color{blue}
\bezier{30}(21,21)(27,21)(33,21)
\bezier{30}(20.6,20.4)(26.6,20.4)(32.6,20.4)
\bezier{30}(20.2,19.8)(26.2,19.8)(32.2,19.8)
\bezier{30}(19.8,19.2)(25.8,19.2)(31.8,19.2)
\bezier{30}(19.4,18.6)(25.4,18.6)(31.4,18.6)
\bezier{30}(19.0,18.0)(25.0,18.0)(31.0,18.0)
\bezier{30}(18.6,17.4)(24.6,17.4)(30.6,17.4)
\bezier{30}(18.2,16.8)(24.2,16.8)(30.2,16.8)
\bezier{30}(17.8,16.2)(23.8,16.2)(29.8,16.2)
\bezier{30}(17.4,15.6)(23.4,15.6)(29.4,15.6)
\bezier{30}(17.0,15.0)(23.0,15.0)(29.0,15.0)
\bezier{30}(16.6,14.4)(22.6,14.4)(28.6,14.4)
\bezier{30}(16.2,13.8)(22.2,13.8)(28.2,13.8)
\bezier{30}(15.8,13.2)(21.8,13.2)(27.8,13.2)
\bezier{30}(15.4,12.6)(21.4,12.6)(27.4,12.6)

\put(20,9.5){\small$\ga_1$}
\put(14.4,17.5){\small$\ga_2$}
\put(15,12){\vector(1,0){12.0}}
\put(15,12){\vector(2,3){6}}
\end{picture}\hspace{8mm}
\begin{picture}(55,53)(0,0)
\put(5,0){\emph{b})}
\put(-2,51){$\bG_\gt$}

\put(21,3){\circle*{1.5}}
\put(33,3){\circle*{1.5}}
\put(45,3){\circle*{1.5}}

\put(27,6){\circle*{1.5}}
\put(39,6){\circle*{1.5}}

\put(15,12){\circle*{1.5}}
\put(27,12){\circle*{1.5}}
\put(39,12){\circle*{1.5}}

\put(21,15){\circle*{1.5}}
\put(33,15){\circle*{1.5}}
\put(45,15){\circle*{1.5}}
\put(21,21){\circle*{1.5}}
\put(33,21){\circle*{1.5}}
\put(45,21){\circle{1.5}}

\put(15,24){\circle*{1.5}}
\put(27,24){\circle*{1.5}}
\put(39,24){\circle*{1.5}}

\put(15,30){\circle*{1.5}}
\put(27,30){\circle*{1.5}}
\put(39,30){\circle*{1.5}}

\put(9,33){\circle*{1.5}}
\put(21,33){\circle*{1.5}}
\put(33,33){\circle*{1.5}}
\put(9,39){\circle*{1.5}}
\put(21,39){\circle*{1.5}}
\put(33,39){\circle*{1.5}}

\put(15,42){\circle*{1.5}}
\put(27,42){\circle*{1.5}}
\put(39,42){\circle{1.5}}

\put(15,48){\circle*{1.5}}
\put(27,48){\circle*{1.5}}

\put(9,51){\circle*{1.5}}
\put(21,51){\circle*{1.5}}
\put(33,51){\circle*{1.5}}
\put(21,3){\line(0,-1){3.00}}
\put(33,3){\line(0,-1){3.00}}
\put(45,3){\line(0,-1){3.00}}

\put(9,51){\line(0,1){3.00}}
\put(21,51){\line(0,1){3.00}}
\put(33,51){\line(0,1){3.00}}

\put(21,3){\line(-2,1){3.90}}
\put(21,3){\line(2,1){6.00}}
\put(33,3){\line(-2,1){6.00}}
\put(33,3){\line(2,1){6.00}}
\put(45,3){\line(-2,1){6.00}}
\put(45,3){\line(2,1){4.4}}

\put(27,6){\line(0,1){6.00}}
\put(39,6){\line(0,1){6.00}}

\put(21,15){\line(-2,-1){6.00}}
\put(21,15){\line(2,-1){6.00}}
\put(33,15){\line(-2,-1){6.00}}
\put(33,15){\line(2,-1){6.00}}
\put(45,15){\line(-2,-1){6.00}}
\bezier{60}(45,15)(46,14.5)(47,14)

\put(21,15){\line(0,1){6.00}}
\put(33,15){\line(0,1){6.00}}
\put(45,15){\line(0,1){6.00}}
\bezier{60}(15,24)(13.55,23.275)(12.1,22.55)
\put(21,21){\line(-2,1){6.00}}
\put(21,21){\line(2,1){6.00}}
\put(33,21){\line(-2,1){6.00}}
\put(33,21){\line(2,1){6.00}}
\put(45,21){\line(-2,1){6.00}}

\put(15,24){\line(0,1){6.00}}
\put(27,24){\line(0,1){6.00}}
\put(39,24){\line(0,1){6.00}}

\put(9,33){\line(2,-1){6.00}}
\put(21,33){\line(-2,-1){6.00}}
\put(21,33){\line(2,-1){6.00}}
\put(33,33){\line(-2,-1){6.00}}
\put(33,33){\line(2,-1){6.00}}
\bezier{60}(39,30)(40.15,30.95)(41.9,31.9)

\put(9,33){\line(0,1){6.00}}
\put(21,33){\line(0,1){6.00}}
\put(33,33){\line(0,1){6.00}}
\bezier{60}(9,39)(8,39.5)(7,40)
\put(9,39){\line(-2,1){2.00}}
\put(9,39){\line(2,1){6.00}}
\put(21,39){\line(-2,1){6.00}}
\put(21,39){\line(2,1){6.00}}
\put(33,39){\line(-2,1){6.00}}
\put(33,39){\line(2,1){6.00}}

\put(15,42){\line(0,1){6.00}}
\put(27,42){\line(0,1){6.00}}
\bezier{60}(9,51)(6.8,49.9)(4.6,48.8)
\put(9,51){\line(2,-1){6.00}}
\put(21,51){\line(-2,-1){6.00}}
\put(21,51){\line(2,-1){6.00}}
\put(33,51){\line(-2,-1){6.00}}
\bezier{60}(33,51)(34.95,50.025)(36.9,49.05)

\put(21.8,15.5){\small$v_1$}
\put(10.5,10.3){\small$v_2$}
\put(33.5,15.5){\small$v_1\!\!+\!\ga_1$}
\put(27.5,10.0){\small$v_2\!\!+\!\ga_1$}

\bezier{100}(3.10,54)(10.81,27)(18.53,0)
\bezier{100}(35.50,54)(43.21,27)(50.93,0)

\color{red}
\bezier{200}(15,12)(30,16.5)(45,21)
\put(44,20.7){\vector(3,1){0.5}}
\put(40,20.3){\small$\gt$}

\color{blue}
\bezier{60}(15.15,11.5)(30.15,16.0)(45.15,20.5)
\bezier{60}(15.30,11.0)(30.30,15.5)(45.30,20.0)
\bezier{60}(15.45,10.5)(30.45,15.0)(45.45,19.5)
\bezier{60}(15.60,10.0)(30.60,14.5)(45.60,19.0)
\bezier{60}(15.75,9.5)(30.75,14.0)(45.75,18.5)
\bezier{60}(15.90,9.0)(30.90,13.5)(45.90,18.0)
\bezier{60}(16.05,8.5)(31.05,13.0)(46.05,17.5)
\put(20,9.5){\small$\ga_1$}
\put(14.4,17.5){\small$\ga_2$}
\put(15,12){\vector(1,0){12.0}}
\put(15,12){\vector(2,3){6}}
\end{picture}
\begin{picture}(27,53)(0,0)
\put(-2,0){\emph{c})}
\put(-3,51){$\bG_\gt$}

\put(9,3){\circle*{1.5}}
\put(21,3){\circle*{1.5}}

\put(3,6){\circle*{1.5}}
\put(15,6){\circle*{1.5}}
\put(27,6){\circle{1.5}}

\put(3,12){\circle*{1.5}}
\put(15,12){\circle*{1.5}}
\put(27,12){\circle{1.5}}

\put(9,15){\circle*{1.5}}
\put(21,15){\circle*{1.5}}
\put(9,21){\circle*{1.5}}
\put(21,21){\circle*{1.5}}

\put(3,24){\circle*{1.5}}
\put(15,24){\circle*{1.5}}
\put(27,24){\circle{1.5}}

\put(3,30){\circle*{1.5}}
\put(15,30){\circle*{1.5}}
\put(27,30){\circle{1.5}}

\put(9,33){\circle*{1.5}}
\put(21,33){\circle*{1.5}}
\put(9,39){\circle*{1.5}}
\put(21,39){\circle*{1.5}}

\put(3,42){\circle*{1.5}}
\put(15,42){\circle*{1.5}}
\put(27,42){\circle{1.5}}

\put(3,48){\circle*{1.5}}
\put(15,48){\circle*{1.5}}
\put(27,48){\circle{1.5}}

\put(9,51){\circle*{1.5}}
\put(21,51){\circle*{1.5}}
\put(9,3){\line(0,-1){3.00}}
\put(21,3){\line(0,-1){3.00}}
\put(9,51){\line(0,1){3.00}}
\put(21,51){\line(0,1){3.00}}

\put(9,3){\line(-2,1){6.00}}
\put(9,3){\line(2,1){6.00}}
\put(21,3){\line(-2,1){6.00}}
\put(21,3){\line(2,1){6.00}}
\put(3,6){\line(0,1){6.00}}
\put(15,6){\line(0,1){6.00}}
\put(27,6){\line(0,1){6.00}}

\put(9,15){\line(-2,-1){6.00}}
\put(9,15){\line(2,-1){6.00}}
\put(21,15){\line(-2,-1){6.00}}
\put(21,15){\line(2,-1){6.00}}

\put(9,15){\line(0,1){6.00}}
\put(21,15){\line(0,1){6.00}}
\put(9,21){\line(-2,1){6.00}}
\put(9,21){\line(2,1){6.00}}
\put(21,21){\line(-2,1){6.00}}
\put(21,21){\line(2,1){6.00}}
\put(3,24){\line(0,1){6.00}}
\put(15,24){\line(0,1){6.00}}
\put(27,24){\line(0,1){6.00}}
\put(9,33){\line(-2,-1){6.00}}
\put(9,33){\line(2,-1){6.00}}
\put(21,33){\line(-2,-1){6.00}}
\put(21,33){\line(2,-1){6.00}}
\put(9,33){\line(0,1){6.00}}
\put(21,33){\line(0,1){6.00}}
\put(9,39){\line(-2,1){6.00}}
\put(9,39){\line(2,1){6.00}}
\put(21,39){\line(-2,1){6.00}}
\put(21,39){\line(2,1){6.00}}

\put(3,42){\line(0,1){6.00}}
\put(15,42){\line(0,1){6.00}}
\put(27,42){\line(0,1){6.00}}

\put(9,51){\line(-2,-1){6.00}}
\put(9,51){\line(2,-1){6.00}}
\put(21,51){\line(-2,-1){6.00}}
\put(21,51){\line(2,-1){6.00}}

\bezier{90}(3,0)(3,27)(3,54)
\bezier{90}(27,0)(27,27)(27,54)
\put(9.8,15.5){\small$v_1$}
\put(-1,10){\small$v_2$}
\put(21.8,15.5){\small$v_3$}
\put(15.2,9.5){\small$v_4$}

\color{red}
\put(23,9.2){\small$\gt$}
\put(3,12){\vector(1,0){24.0}}

\color{blue}
\put(8,9.5){\small$\ga_1$}
\put(2.9,17.9){\small$\ga_2$}
\put(3,11.8){\vector(1,0){12.0}}
\put(3,12){\vector(2,3){6}}
\bezier{50}(3,12.5)(15,12.5)(27,12.5)
\bezier{50}(3,13.0)(15,13.0)(27,13.0)
\bezier{50}(3,13.5)(15,13.5)(27,13.5)
\bezier{50}(3,14.0)(15,14.0)(27,14.0)
\bezier{50}(3,14.5)(15,14.5)(27,14.5)
\bezier{50}(3,15.0)(15,15.0)(27,15.0)
\bezier{50}(3,15.5)(15,15.5)(27,15.5)
\bezier{50}(3,16.0)(15,16.0)(27,16.0)
\bezier{50}(3,16.5)(15,16.5)(27,16.5)
\bezier{50}(3,17.0)(15,17.0)(27,17.0)
\bezier{50}(3,17.5)(15,17.5)(27,17.5)
\bezier{50}(3,18.0)(15,18.0)(27,18.0)
\bezier{50}(3,18.5)(15,18.5)(27,18.5)
\bezier{50}(3,19.0)(15,19.0)(27,19.0)
\bezier{50}(3,19.5)(15,19.5)(27,19.5)
\bezier{50}(3,20.0)(15,20.0)(27,20.0)
\bezier{50}(3,20.5)(15,20.5)(27,20.5)
\bezier{50}(3,21.0)(15,21.0)(27,21.0)
\end{picture}

\begin{picture}(23,25)(0,0)
\put(0,0){\emph{d})}
\put(10,5){\circle*{1.5}}
\put(28,14){\circle*{1.5}}
\put(10,5){\line(2,1){18.00}}
\put(6.7,2.5){$v_2$}
\put(28.8,14){$v_1$}
\bezier{200}(10,5)(25,-2.0)(28,14.0)
\bezier{200}(10,5)(12,22)(28,14)
\put(6,16){$\bG_*$}
\end{picture}\hspace{30mm}
\begin{picture}(46,25)(0,0)
\put(0,0){\emph{e})}
\put(10,5){\circle*{1.5}}
\put(22,11){\circle*{1.5}}
\put(10,5){\line(2,1){12.00}}
\put(6.7,2.5){$v_2$}
\put(22,12){$v_1$}
\put(22,11){\line(2,-1){12.00}}
\put(34,5){\circle*{1.5}}
\put(46,11){\circle*{1.5}}
\put(34,5){\line(2,1){12.00}}
\put(30.7,2.5){$v_4$}
\put(46,12){$v_3$}
\bezier{300}(10,5)(36,-6.0)(46,11.0)
\bezier{200}(10,5)(12,15)(22,11)
\bezier{200}(34,5)(36,15)(46,11)
\put(6,16){$\wt\bG_*$}
\end{picture}

\caption{\scriptsize\emph{a}) The hexagonal lattice $\bG$; the vectors $\ga_1,\ga_2$ (marked in blue color) are the periods of $\bG$; the fundamental cell $\Omega$ is spanned by the vectors $\ga_1$ and $\ga_2$. \emph{b})~The nanotube $\bG_\gt$ with the red chiral vector $\gt=2\ga_1+\ga_2$. The dashed lines have to be identified. \emph{c})~The zigzag nanotube $\bG_\gt$ with $\gt=2\ga_1$. \emph{d})~The fundamental graph $\bG_*$ of the hexagonal lattice $\bG$ and any nanotube $\bG_\gt$ with a primitive chiral vector~$\gt$. \emph{e})~The fundamental graph $\wt\bG_*$ of the nanotube $\bG_\gt$ with the non-primitive vector $\gt=2\ga_1$.}
\label{fig1}
\end{figure}

\begin{example}[\textbf{Nanotubes}] The hexagonal lattice $\bG$ is a $\Z^2$-periodic graph, the vectors $\ga_1,\ga_2$ are the fixed periods of $\bG$, see Fig.~\ref{fig1}\emph{a}. For non-zero vectors $\gt\in\Z^2$ the subcoverings $\bG_\gt=\bG/\Z\gt$ of the hexagonal lattice $\bG$ are known as \emph{nanotubes}. For example, the nanotube $\bG_\gt$ with the chiral vector $\gt=(2,1)$ is shown in Fig. \ref{fig1}\emph{b}. The so-called zigzag nanotube $\bG_\gt$, $\gt=(2,0)$, is presented in Fig. \ref{fig1}\emph{c}.
\end{example}

\begin{figure}[t!]\centering
\unitlength 1.0mm 
\linethickness{0.4pt}

\begin{picture}(46,50)(0,0)
\put(2,8){\emph{a})}
\put(7.9,17.5){$\scriptstyle v$}
\put(5,42){$\dL^3$}
\put(8,10){\line(1,0){34.00}}
\put(8,20){\line(1,0){34.00}}
\put(8,30){\line(1,0){34.00}}
\put(8,40){\line(1,0){34.00}}
\put(10,8){\line(0,1){34.00}}
\put(20,8){\line(0,1){34.00}}
\put(30,8){\line(0,1){34.00}}
\put(40,8){\line(0,1){34.00}}

\put(14,14){\line(1,0){34.00}}
\put(14,24){\line(1,0){34.00}}
\put(14,34){\line(1,0){34.00}}
\put(14,44){\line(1,0){34.00}}
\put(16,12){\line(0,1){34.00}}
\put(26,12){\line(0,1){34.00}}
\put(36,12){\line(0,1){34.00}}
\put(46,12){\line(0,1){34.00}}

\put(20.5,18){\line(1,0){33.5}}
\put(20.5,28){\line(1,0){33.5}}
\put(20.5,38){\line(1,0){33.5}}
\put(20,48){\line(1,0){34.00}}
\put(22,16){\line(0,1){34.00}}
\put(32,16){\line(0,1){34.00}}
\put(42,16){\line(0,1){34.00}}
\put(52,16){\line(0,1){34.00}}

\put(16.5,26.5){$\scriptstyle\Omega$}

\put(10,10){\circle*{1.5}}
\put(20,10){\circle*{1.5}}
\put(30,10){\circle*{1.5}}
\put(40,10){\circle*{1.5}}
\put(10,20){\circle*{1.5}}
\put(20,20){\circle*{1.5}}
\put(30,20){\circle*{1.5}}
\put(40,20){\circle*{1.5}}
\put(10,30){\circle*{1.5}}
\put(20,30){\circle*{1.5}}
\put(30,30){\circle*{1.5}}
\put(40,30){\circle*{1.5}}
\put(10,40){\circle*{1.5}}
\put(20,40){\circle*{1.5}}
\put(30,40){\circle*{1.5}}
\put(40,40){\circle*{1.5}}

\put(16,14){\circle*{1.5}}
\put(26,14){\circle*{1.5}}
\put(36,14){\circle*{1.5}}
\put(46,14){\circle*{1.5}}
\put(16,24){\circle*{1.5}}
\put(26,24){\circle*{1.5}}
\put(36,24){\circle*{1.5}}
\put(46,24){\circle*{1.5}}
\put(16,34){\circle*{1.5}}
\put(26,34){\circle*{1.5}}
\put(36,34){\circle*{1.5}}
\put(46,34){\circle*{1.5}}
\put(16,44){\circle*{1.5}}
\put(26,44){\circle*{1.5}}
\put(36,44){\circle*{1.5}}
\put(46,44){\circle*{1.5}}

\put(22,18){\circle*{1.5}}
\put(32,18){\circle*{1.5}}
\put(42,18){\circle*{1.5}}
\put(52,18){\circle*{1.5}}
\put(22,28){\circle*{1.5}}
\put(32,28){\circle*{1.5}}
\put(42,28){\circle*{1.5}}
\put(52,28){\circle*{1.5}}
\put(22,38){\circle*{1.5}}
\put(32,38){\circle*{1.5}}
\put(42,38){\circle*{1.5}}
\put(52,38){\circle*{1.5}}
\put(22,48){\circle*{1.5}}
\put(32,48){\circle*{1.5}}
\put(42,48){\circle*{1.5}}
\put(52,48){\circle*{1.5}}

\put(8.2,8.8){\line(3,2){15.6}}
\put(18.2,8.8){\line(3,2){15.6}}
\put(28.2,8.8){\line(3,2){15.6}}
\put(38.2,8.8){\line(3,2){15.6}}
\put(8.2,18.8){\line(3,2){15.6}}
\put(18.2,18.8){\line(3,2){15.6}}
\put(28.2,18.8){\line(3,2){15.6}}
\put(38.2,18.8){\line(3,2){15.6}}
\put(8.2,28.8){\line(3,2){15.6}}
\put(18.2,28.8){\line(3,2){15.6}}
\put(28.2,28.8){\line(3,2){15.6}}
\put(38.2,28.8){\line(3,2){15.6}}
\put(8.2,38.8){\line(3,2){15.6}}
\put(18.2,38.8){\line(3,2){15.6}}
\put(28.2,38.8){\line(3,2){15.6}}
\put(38.2,38.8){\line(3,2){15.6}}
\color{red}
\put(25,14.5){\vector(3,-1){0.5}}
\put(33,25){$\scriptstyle \gt_2$}
\put(23.0,15.7){$\scriptstyle \gt_1$}
\put(35,24){\vector(3,1){0.5}}
\bezier{300}(10,20)(18,17)(26,14)
\bezier{300}(10,20)(23,22)(36,24)
\color{blue}
\put(10,20){\vector(1,0){10.00}}
\put(10,20){\vector(0,1){10.00}}
\put(10,20){\vector(3,2){6.00}}
\put(16,18.3){$\scriptstyle \ga_1$}
\put(11,23){$\scriptstyle  \ga_2$}
\put(6.7,27){$\scriptstyle \ga_3$}
\bezier{15}(10.5,20)(10.5,25)(10.5,30)
\bezier{15}(11,20)(11,25)(11,30)
\bezier{15}(11.5,20)(11.5,25)(11.5,30)
\bezier{15}(12,20)(12,25)(12,30)
\bezier{15}(12.5,20)(12.5,25)(12.5,30)
\bezier{15}(13,20)(13,25)(13,30)
\bezier{15}(13.5,20)(13.5,25)(13.5,30)
\bezier{15}(14,20)(14,25)(14,30)
\bezier{15}(14.5,20)(14.5,25)(14.5,30)
\bezier{15}(15,20)(15,25)(15,30)
\bezier{15}(15.5,20)(15.5,25)(15.5,30)
\bezier{15}(16,20)(16,25)(16,30)
\bezier{15}(16.5,20)(16.5,25)(16.5,30)
\bezier{15}(17,20)(17,25)(17,30)
\bezier{15}(17.5,20)(17.5,25)(17.5,30)
\bezier{15}(18,20)(18,25)(18,30)
\bezier{15}(18.5,20)(18.5,25)(18.5,30)
\bezier{15}(19,20)(19,25)(19,30)
\bezier{15}(19.5,20)(19.5,25)(19.5,30)

\bezier{15}(20.6,20.4)(20.6,25.4)(20.6,30.4)
\bezier{15}(21.2,20.8)(21.2,25.8)(21.2,30.8)
\bezier{15}(21.8,21.2)(21.8,26.2)(21.8,31.2)
\bezier{15}(22.4,21.6)(22.4,26.6)(22.4,31.6)
\bezier{15}(23.0,22.0)(23.0,27.0)(23.0,32.0)
\bezier{15}(23.6,22.4)(23.6,27.4)(23.6,32.4)
\bezier{15}(24.2,22.8)(24.2,27.8)(24.2,32.8)
\bezier{15}(24.8,23.2)(24.8,28.2)(24.8,33.2)
\bezier{15}(25.4,23.6)(25.4,28.6)(25.4,33.6)
\bezier{15}(26.0,24.0)(26.0,29.0)(26.0,34.0)

\bezier{15}(10.6,30.4)(15.6,30.4)(20.6,30.4)
\bezier{15}(11.2,30.8)(16.2,30.8)(21.2,30.8)
\bezier{15}(11.8,31.2)(16.8,31.2)(21.8,31.2)
\bezier{15}(12.4,31.6)(17.4,31.6)(22.4,31.6)
\bezier{15}(13.0,32.0)(18.0,32.0)(23.0,32.0)
\bezier{15}(13.6,32.4)(18.6,32.4)(23.6,32.4)
\bezier{15}(14.2,32.8)(19.2,32.8)(24.2,32.8)
\bezier{15}(14.8,33.2)(19.8,33.2)(24.8,33.2)
\bezier{15}(15.4,33.6)(20.4,33.6)(25.4,33.6)
\end{picture}\hspace{10mm}
\begin{picture}(42,40)
\put(7,10){\line(1,0){36.00}}
\put(7,20){\line(1,0){36.00}}
\put(7,30){\line(1,0){36.00}}
\put(7,40){\line(1,0){36.00}}
\put(10,7){\line(0,1){36.00}}
\put(20,7){\line(0,1){36.00}}
\put(30,7){\line(0,1){36.00}}
\put(40,7){\line(0,1){36.00}}
\put(10.2,8){$\scriptstyle v$}

\put(7,7){\line(1,1){36.00}}
\put(7,17){\line(1,1){26.00}}
\put(17,7){\line(1,1){26.00}}
\put(7,27){\line(1,1){16.00}}
\put(27,7){\line(1,1){16.00}}
\put(37,7){\line(1,1){6.00}}
\put(7,37){\line(1,1){6.00}}

\put(10,10){\circle*{1.5}}
\put(20,10){\circle*{1.5}}
\put(30,10){\circle*{1.5}}
\put(40,10){\circle*{1.5}}

\put(10,20){\circle*{1.5}}
\put(20,20){\circle*{1.5}}
\put(30,20){\circle*{1.5}}
\put(40,20){\circle*{1.5}}

\put(10,30){\circle*{1.5}}
\put(20,30){\circle*{1.5}}
\put(30,30){\circle*{1.5}}
\put(40,30){\circle*{1.5}}

\put(10,40){\circle*{1.5}}
\put(20,40){\circle*{1.5}}
\put(30,40){\circle*{1.5}}
\put(40,40){\circle*{1.5}}
\put(20,20){\line(1,1){10.00}}
\put(2,8){\emph{b})}
\put(4,42){$\dL^3_{\gt_1}$}
\color{red}
\put(10,10){\vector(2,1){20.00}}
\put(22,17.6){$\scriptstyle \gt_2$}
\color{blue}
\put(10,10){\vector(1,0){10.00}}
\put(10,10){\vector(0,1){10.00}}
\put(14,8){$\scriptstyle\ga_1$}
\put(6.5,15){$\scriptstyle \ga_2$}
\bezier{15}(10.5,10)(10.5,15)(10.5,20)
\bezier{15}(11,10)(11,15)(11,20)
\bezier{15}(11.5,10)(11.5,15)(11.5,20)
\bezier{15}(12,10)(12,15)(12,20)
\bezier{15}(12.5,10)(12.5,15)(12.5,20)
\bezier{15}(13,10)(13,15)(13,20)
\bezier{15}(13.5,10)(13.5,15)(13.5,20)
\bezier{15}(14,10)(14,15)(14,20)
\bezier{15}(14.5,10)(14.5,15)(14.5,20)
\bezier{15}(15,10)(15,15)(15,20)
\bezier{15}(15.5,10)(15.5,15)(15.5,20)
\bezier{15}(16,10)(16,15)(16,20)
\bezier{15}(16.5,10)(16.5,15)(16.5,20)
\bezier{15}(17,10)(17,15)(17,20)
\bezier{15}(17.5,10)(17.5,15)(17.5,20)
\bezier{15}(18,10)(18,15)(18,20)
\bezier{15}(18.5,10)(18.5,15)(18.5,20)
\bezier{15}(19,10)(19,15)(19,20)
\bezier{15}(19.5,10)(19.5,15)(19.5,20)
\end{picture}\;
\hspace{5mm}
\begin{picture}(30,45)
\put(0,8){\emph{d})}
\put(6,28){$\dL^3_*$}
\put(17,11.5){$v$}
\put(18.0,15){\circle*{1.5}}
\bezier{200}(18.0,15)(5,22)(4.0,15)
\bezier{200}(18.0,15)(5,8)(4.0,15)
\bezier{200}(18.0,15)(31,22)(32.0,15)
\bezier{200}(18.0,15)(31,9)(32.0,15)
\bezier{200}(18.0,15)(25,28)(18.0,29)
\bezier{200}(18.0,15)(11,28)(18.0,29)
\end{picture}

\begin{picture}(132,15)
\put(5,10){\line(1,0){100.00}} \put(10,10){\circle*{1.5}}
\put(20,10){\circle*{1.5}} \put(30,10){\circle*{1.5}}
\put(40,10){\circle*{1.5}} \put(50,10){\circle*{1.5}}
\put(60,10){\circle*{1.5}} \put(70,10){\circle*{1.5}} \put(80,10){\circle*{1.5}}
\put(90,10){\circle*{1.5}} \put(100,10){\circle*{1.5}}
\put(9,6.5){$\scriptstyle 0$}
\put(19,6.5){$\scriptstyle 1$} \put(29,6.5){$\scriptstyle 2$} \put(39,6.5){$\scriptstyle3$}
\put(49,6.5){$\scriptstyle 4$} \put(59,6.5){$\scriptstyle 5$} \put(69,6.5){$\scriptstyle 6$}\put(79,6.5){$\scriptstyle 7$}
\put(89,6.5){$\scriptstyle 8$}\put(99,6.5){$\scriptstyle 9$}

\bezier{600}(5,13.8)(7,13)(10,10)
\bezier{600}(5,13.8)(13,17.5)(20,10)
\bezier{600}(10,10)(20,20)(30,10)
\bezier{600}(20,10)(30,20)(40,10)
\bezier{600}(30,10)(40,20)(50,10)
\bezier{600}(40,10)(50,20)(60,10)
\bezier{600}(50,10)(60,20)(70,10)
\bezier{600}(60,10)(70,20)(80,10)
\bezier{600}(70,10)(80,20)(90,10)
\bezier{600}(80,10)(90,20)(100,10)
\bezier{600}(90,10)(97,17.5)(105,13.8)
\bezier{600}(100,10)(103,13)(105,13.8)

\bezier{600}(5,6.5)(8.3,6.9)(10,10)
\bezier{600}(10,10)(15,3)(20,10)
\bezier{600}(20,10)(25,3)(30,10)
\bezier{600}(30,10)(35,3)(40,10)
\bezier{600}(40,10)(45,3)(50,10)
\bezier{600}(50,10)(55,3)(60,10)
\bezier{600}(60,10)(65,3)(70,10)
\bezier{600}(70,10)(75,3)(80,10)
\bezier{600}(80,10)(85,3)(90,10)
\bezier{600}(90,10)(95,3)(100,10)
\bezier{600}(100,10)(101.7,6.9)(105,6.5)

\put(0,12){$\dL^3_\gt$}
\put(0,6){\emph{c})}
\color{blue}
\put(10,10){\vector(1,0){10}}
\put(14,8){$\scriptstyle\ga_1$}
\end{picture}

\caption{\scriptsize\emph{a}) The cubic lattice $\dL^3$; the blue vectors $\ga_1,\ga_2,\ga_3$ are the periods of $\dL^3$; the fundamental cell $\Omega$ is shaded. The chiral vectors $\gt_1=\ga_1+\ga_2-\ga_3$ and $\gt_2=2\ga_1+\ga_2$ are marked in red color. \emph{b}) The subcovering $\dL^3_{\gt_1}$ of $\dL^3$ with the chiral vector $\gt_1$ (the triangular lattice). \emph{c}) The subcovering $\dL^3_{\gt}$ of $\dL^3$ with $\gt=\{\gt_1,\gt_2\}$. \emph{d}) The fundamental graph $\dL^3_*$ of the cubic lattice $\dL^3$ and any its subcovering $\dL^3_\gt$ with a primitive set $\gt$ of chiral vectors.} \label{fig2}
\end{figure}

\begin{example} \lb{Excl} The cubic lattice $\dL^3$ is a $\Z^3$-periodic graph, the vectors $\ga_1,\ga_2,\ga_3$ are the periods of $\dL^3$, see Fig.~\ref{fig2}\emph{a}. The subcovering $\dL^3_{\gt_1}=\dL^3/\Z\gt_1$ of the cubic lattice $\dL^3$ with the chiral vector $\gt_1=(1,1,-1)$ is the triangular lattice, Fig. \ref{fig2}\emph{b}. Let $\G_\gt$ be the sublattice of $\Z^3$ with the basis $\gt:=\{\gt_1,\gt_2\}$, where $\gt_1=(1,1,-1)$ and $\gt_2=(2,1,0)$. Then the subcovering $\dL^3_{\gt}=\dL^3/\G_\gt$ of the cubic lattice $\dL^3$ is the $\Z$-periodic graph with the period $\ga_1$ shown in Fig. \ref{fig2}\emph{c}. Note that this $\Z$-periodic graph can also be considered as the subcovering of the triangular lattice (which is a $\Z^2$-periodic graph with the periods $\ga_1,\ga_2$) with the chiral vector $\gt_2=(2,1)$.
\end{example}

\begin{figure}[h]
\unitlength 1.0mm 
\linethickness{0.4pt}
\unitlength 0.9mm
\begin{picture}(70,60)
\put(0,-2){\emph{a})}
\put(5,53){$\cG$}

\put(20,20.5){\tiny$v_1$}
\put(8.8,21.6){\tiny$v_2$}
\put(8.1,41.9){\tiny$v_2\!+\!\ga_3$}

\put(60,30){\tiny$\bG$}
\put(57,10){\tiny$\bG\!-\!\ga_3$}
\put(57,50){\tiny$\bG\!+\!\ga_3$}

\qbezier(21,3)(16.5,11.5)(12,20)
\qbezier(39,3)(34.5,11.5)(30,20)
\qbezier(57,3)(52.5,11.5)(48,20)
\qbezier(15,9)(10.5,17.5)(6,26)
\qbezier(33,9)(28.5,17.5)(24,26)
\qbezier(51,9)(46.5,17.5)(42,26)
\qbezier(69,9)(64.5,17.5)(60,26)
\qbezier(27,15)(22.5,23.5)(18,32)
\qbezier(45,15)(40.5,23.5)(36,32)
\qbezier(63,15)(58.5,23.5)(54,32)

\qbezier(21,23)(16.5,31.5)(12,40)
\qbezier(39,23)(34.5,31.5)(30,40)
\qbezier(57,23)(52.5,31.5)(48,40)
\qbezier(15,29)(10.5,37.5)(6,46)
\qbezier(33,29)(28.5,37.5)(24,46)
\qbezier(51,29)(46.5,37.5)(42,46)
\qbezier(69,29)(64.5,37.5)(60,46)
\qbezier(27,35)(22.5,43.5)(18,52)
\qbezier(45,35)(40.5,43.5)(36,52)
\qbezier(63,35)(58.5,43.5)(54,52)

\qbezier(12,0)(12.5,-1)(13,-2)
\qbezier(30,0)(30.5,-1)(31,-2)
\qbezier(48,0)(48.5,-1)(49,-2)
\qbezier(6,6)(6.5,5)(7,4)
\qbezier(24,6)(24.5,5)(25,4)
\qbezier(42,6)(42.5,5)(43,4)
\qbezier(60,6)(60.5,5)(61,4)
\qbezier(18,12)(18.5,11)(19,10)
\qbezier(36,12)(36.5,11)(37,10)
\qbezier(54,12)(54.5,11)(55,10)
\qbezier(72,12)(72.5,11)(73,10)
\qbezier(72,32)(72.5,31)(73,30)
\qbezier(72,52)(72.5,51)(73,50)

\qbezier(3,3)(2.5,4)(2,5)
\qbezier(3,23)(2.5,24)(2,25)
\qbezier(3,43)(2.5,44)(2,45)
\qbezier(21,43)(20.5,44)(20,45)
\qbezier(39,43)(38.5,44)(38,45)
\qbezier(57,43)(56.5,44)(56,45)
\qbezier(15,49)(14.5,50)(14,51)
\qbezier(33,49)(32.5,50)(32,51)
\qbezier(51,49)(50.5,50)(50,51)
\qbezier(69,49)(68.5,50)(68,51)

\qbezier(27,55)(26.5,56)(26,57)
\qbezier(45,55)(44.5,56)(44,57)
\qbezier(63,55)(62.5,56)(62,57)

\qbezier(12,0)(11.5,-0.5)(11,-1)
\qbezier(30,0)(29.5,-0.5)(29,-1)
\qbezier(48,0)(47.5,-0.5)(47,-1)
\qbezier(3,3)(2.25,2.75)(1.5,2.5)
\qbezier(57,3)(57.75,2.75)(58.5,2.5)
\qbezier(6,6)(5.25,6.25)(4.5,6.5)
\qbezier(69,9)(69.75,8.75)(70.5,8.5)
\qbezier(18,12)(17.25,12.25)(16.5,12.5)
\qbezier(72,12)(72.75,12.25)(73.5,12.5)
\qbezier(27,15)(27.5,15.5)(28,16)
\qbezier(45,15)(45.5,15.5)(46,16)
\qbezier(63,15)(63.5,15.5)(64,16)

\qbezier(12,20)(11.5,19.5)(11,19)
\qbezier(30,20)(29.5,19.5)(29,19)
\qbezier(48,20)(47.5,19.5)(47,19)
\qbezier(3,23)(2.25,22.75)(1.5,22.5)
\qbezier(57,23)(57.75,22.75)(58.5,22.5)
\qbezier(6,26)(5.25,26.25)(4.5,26.5)
\qbezier(69,29)(69.75,28.75)(70.5,28.5)
\qbezier(18,32)(17.25,32.25)(16.5,32.5)
\qbezier(72,32)(72.75,32.25)(73.5,32.5)
\qbezier(27,35)(27.5,35.5)(28,36)
\qbezier(45,35)(45.5,35.5)(46,36)
\qbezier(63,35)(63.5,35.5)(64,36)

\qbezier(12,40)(11.5,39.5)(11,39)
\qbezier(30,40)(29.5,39.5)(29,39)
\qbezier(48,40)(47.5,39.5)(47,39)
\qbezier(3,43)(2.25,42.75)(1.5,42.5)
\qbezier(57,43)(57.75,42.75)(58.5,42.5)
\qbezier(6,46)(5.25,46.25)(4.5,46.5)
\qbezier(69,49)(69.75,48.75)(70.5,48.5)
\qbezier(18,52)(17.25,52.25)(16.5,52.5)
\qbezier(72,52)(72.75,52.25)(73.5,52.5)
\qbezier(27,55)(27.5,55.5)(28,56)
\qbezier(45,55)(45.5,55.5)(46,56)
\qbezier(63,55)(63.5,55.5)(64,56)

\put(12,0){\circle*{1.5}}
\put(30,0){\circle*{1.5}}
\put(48,0){\circle*{1.5}}
\put(3,3){\circle{1.5}}
\put(21,3){\circle{1.5}}
\put(39,3){\circle{1.5}}
\put(57,3){\circle{1.5}}
\put(6,6){\circle*{1.5}}
\put(24,6){\circle*{1.5}}
\put(42,6){\circle*{1.5}}
\put(60,6){\circle*{1.5}}
\put(15,9){\circle{1.5}}
\put(33,9){\circle{1.5}}
\put(51,9){\circle{1.5}}
\put(69,9){\circle{1.5}}
\put(18,12){\circle*{1.5}}
\put(36,12){\circle*{1.5}}
\put(54,12){\circle*{1.5}}
\put(72,12){\circle*{1.5}}
\put(27,15){\circle{1.5}}
\put(45,15){\circle{1.5}}
\put(63,15){\circle{1.5}}

\put(12,0){\line(-3,1){9.00}}
\put(12,0){\line(3,1){9.00}}
\put(30,0){\line(-3,1){9.00}}
\put(30,0){\line(3,1){9.00}}
\put(48,0){\line(-3,1){9.00}}
\put(48,0){\line(3,1){9.00}}

\qbezier(6,6)(4.5,4.5)(3,3)
\put(6,6){\line(3,1){9.00}}
\qbezier(24,6)(22.5,4.5)(21,3)
\put(24,6){\line(3,1){9.00}}
\put(24,6){\line(-3,1){9.00}}
\qbezier(42,6)(40.5,4.5)(39,3)
\put(42,6){\line(3,1){9.00}}
\put(42,6){\line(-3,1){9.00}}
\qbezier(60,6)(58.5,4.5)(57,3)
\put(60,6){\line(3,1){9.00}}
\put(60,6){\line(-3,1){9.00}}
\qbezier(18,12)(16.5,10.5)(15,9)
\put(18,12){\line(3,1){9.00}}
\qbezier(36,12)(34.5,10.5)(33,9)
\put(36,12){\line(3,1){9.00}}
\put(36,12){\line(-3,1){9.00}}
\qbezier(54,12)(52.5,10.5)(51,9)
\put(54,12){\line(3,1){9.00}}
\put(54,12){\line(-3,1){9.00}}
\qbezier(72,12)(70.5,10.5)(69,9)
\put(72,12){\line(-3,1){9.00}}


\put(12,20){\circle*{1.5}}
\put(30,20){\circle*{1.5}}
\put(48,20){\circle*{1.5}}
\put(3,23){\circle{1.5}}
\put(21,23){\circle{1.5}}
\put(39,23){\circle{1.5}}
\put(57,23){\circle{1.5}}
\put(6,26){\circle*{1.5}}
\put(24,26){\circle*{1.5}}
\put(42,26){\circle*{1.5}}
\put(60,26){\circle*{1.5}}
\put(15,29){\circle{1.5}}
\put(33,29){\circle{1.5}}
\put(51,29){\circle{1.5}}
\put(69,29){\circle{1.5}}
\put(18,32){\circle*{1.5}}
\put(36,32){\circle*{1.5}}
\put(54,32){\circle*{1.5}}
\put(72,32){\circle*{1.5}}
\put(27,35){\circle{1.5}}
\put(45,35){\circle{1.5}}
\put(63,35){\circle{1.5}}

\put(12,20){\line(-3,1){9.00}}
\put(12,20){\line(3,1){9.00}}
\put(30,20){\line(-3,1){9.00}}
\put(30,20){\line(3,1){9.00}}
\put(48,20){\line(-3,1){9.00}}
\put(48,20){\line(3,1){9.00}}

\qbezier(6,26)(4.5,24.5)(3,23)
\put(6,26){\line(3,1){9.00}}
\qbezier(24,26)(22.5,24.5)(21,23)
\put(24,26){\line(3,1){9.00}}
\put(24,26){\line(-3,1){9.00}}
\qbezier(42,26)(40.5,24.5)(39,23)
\put(42,26){\line(3,1){9.00}}
\put(42,26){\line(-3,1){9.00}}
\qbezier(60,26)(58.5,24.5)(57,23)
\put(60,26){\line(3,1){9.00}}
\put(60,26){\line(-3,1){9.00}}
\qbezier(18,32)(16.5,30.5)(15,29)
\put(18,32){\line(3,1){9.00}}
\qbezier(36,32)(34.5,30.5)(33,29)
\put(36,32){\line(3,1){9.00}}
\put(36,32){\line(-3,1){9.00}}
\qbezier(54,32)(52.5,30.5)(51,29)
\put(54,32){\line(3,1){9.00}}
\put(54,32){\line(-3,1){9.00}}
\qbezier(72,32)(70.5,30.5)(69,29)
\put(72,32){\line(-3,1){9.00}}


\put(12,40){\circle*{1.5}}
\put(30,40){\circle*{1.5}}
\put(48,40){\circle*{1.5}}
\put(3,43){\circle{1.5}}
\put(21,43){\circle{1.5}}
\put(39,43){\circle{1.5}}
\put(57,43){\circle{1.5}}
\put(6,46){\circle*{1.5}}
\put(24,46){\circle*{1.5}}
\put(42,46){\circle*{1.5}}
\put(60,46){\circle*{1.5}}
\put(15,49){\circle{1.5}}
\put(33,49){\circle{1.5}}
\put(51,49){\circle{1.5}}
\put(69,49){\circle{1.5}}
\put(18,52){\circle*{1.5}}
\put(36,52){\circle*{1.5}}
\put(54,52){\circle*{1.5}}
\put(72,52){\circle*{1.5}}
\put(27,55){\circle{1.5}}
\put(45,55){\circle{1.5}}
\put(63,55){\circle{1.5}}

\put(12,40){\line(-3,1){9.00}}
\put(12,40){\line(3,1){9.00}}
\put(30,40){\line(-3,1){9.00}}
\put(30,40){\line(3,1){9.00}}
\put(48,40){\line(-3,1){9.00}}
\put(48,40){\line(3,1){9.00}}

\qbezier(6,46)(4.5,44.5)(3,43)
\put(6,46){\line(3,1){9.00}}
\qbezier(24,46)(22.5,44.5)(21,43)
\put(24,46){\line(3,1){9.00}}
\put(24,46){\line(-3,1){9.00}}
\qbezier(42,46)(40.5,44.5)(39,43)
\put(42,46){\line(3,1){9.00}}
\put(42,46){\line(-3,1){9.00}}
\qbezier(60,46)(58.5,44.5)(57,43)
\put(60,46){\line(3,1){9.00}}
\put(60,46){\line(-3,1){9.00}}

\qbezier(18,52)(16.5,50.5)(15,49)
\put(18,52){\line(3,1){9.00}}
\qbezier(36,52)(34.5,50.5)(33,49)
\put(36,52){\line(3,1){9.00}}
\put(36,52){\line(-3,1){9.00}}
\qbezier(54,52)(52.5,50.5)(51,49)
\put(54,52){\line(3,1){9.00}}
\put(54,52){\line(-3,1){9.00}}
\qbezier(72,52)(70.5,50.5)(69,49)
\put(72,52){\line(-3,1){9.00}}

\put(30,23){\tiny$\Omega$}
\color{red}
\qbezier(12,20)(27,13)(42,6)
\put(41,6.5){\vector(2,-1){1.00}}
\put(39,8){\tiny$\gt$}

\color{blue}
\put(12,20){\vector(2,1){12.00}}
\put(12,20){\vector(1,0){18.00}}
\put(12,20){\vector(0,1){20.00}}

\bezier{30}(13.0,20.5)(22.0,20.5)(31.0,20.5)
\bezier{30}(14.0,21.0)(23.0,21.0)(32.0,21.0)
\bezier{30}(15.0,21.5)(24.0,21.5)(33.0,21.5)
\bezier{30}(16.0,22.0)(25.0,22.0)(34.0,22.0)
\bezier{30}(17.0,22.5)(26.0,22.5)(35.0,22.5)
\bezier{30}(18.0,23.0)(27.0,23.0)(36.0,23.0)
\bezier{30}(19.0,23.5)(28.0,23.5)(37.0,23.5)
\bezier{30}(20.0,24.0)(29.0,24.0)(38.0,24.0)
\bezier{30}(21.0,24.5)(30.0,24.5)(39.0,24.5)
\bezier{30}(22.0,25.0)(31.0,25.0)(40.0,25.0)
\bezier{30}(23.0,25.5)(32.0,25.5)(41.0,25.5)
\bezier{30}(24.0,26.0)(33.0,26.0)(42.0,26.0)

\put(20,18){\tiny$\ga_1$}
\put(15,24){\tiny$\ga_2$}
\put(8.5,33){\tiny$\ga_3$}
\color{green}
\qbezier(12,20)(21.2,10.2)(30.4,0.4)
\put(29.6,1.1){\vector(1,-1){1.00}}
\put(29.5,1.7){\tiny$\gt_1$}
\qbezier(12,20)(18,13)(24,6)
\put(22.7,7.2){\vector(1,-1){1.00}}
\put(19.0,8){\tiny$\gt_2$}
\end{picture}
\unitlength 1.0mm
\hspace{10mm}\begin{picture}(40,45)(0,0)
\put(-4,0){\emph{b})}
\put(-2,38){$\cG_\gt$}

\put(9,3){\circle*{1.5}}
\put(21,3){\circle*{1.5}}
\put(33,3){\circle*{1.5}}

\put(3,6){\circle{1.5}}
\put(15,6){\circle{1.5}}
\put(27,6){\circle{1.5}}
\put(39,6){\circle{1.5}}

\put(3,12){\circle*{1.5}}
\put(15,12){\circle*{1.5}}
\put(27,12){\circle*{1.5}}
\put(39,12){\circle*{1.5}}

\put(9,15){\circle{1.5}}
\put(21,15){\circle{1.5}}
\put(33,15){\circle{1.5}}
\put(9,21){\circle*{1.5}}
\put(21,21){\circle*{1.5}}
\put(33,21){\circle*{1.5}}

\put(3,24){\circle{1.5}}
\put(15,24){\circle{1.5}}
\put(27,24){\circle{1.5}}
\put(39,24){\circle{1.5}}

\put(3,30){\circle*{1.5}}
\put(15,30){\circle*{1.5}}
\put(27,30){\circle*{1.5}}
\put(39,30){\circle*{1.5}}

\put(9,33){\circle{1.5}}
\put(21,33){\circle{1.5}}
\put(33,33){\circle{1.5}}
\put(9,3){\line(0,-1){3.00}}
\put(21,3){\line(0,-1){3.00}}
\put(33,3){\line(0,-1){3.00}}

\put(9,33){\line(0,1){3.00}}
\put(21,33){\line(0,1){3.00}}
\put(33,33){\line(0,1){3.00}}

\bezier{60}(3,6)(1.5,5.25)(0,4.5)
\bezier{60}(3,12)(1.5,12.75)(0,13.5)
\bezier{60}(3,24)(1.5,23.25)(0,22.5)
\bezier{60}(3,30)(1.5,30.75)(0,31.5)

\bezier{60}(39,6)(40.5,5.25)(42,4.5)
\bezier{60}(39,12)(40.5,12.75)(42,13.5)
\bezier{60}(39,24)(40.5,23.25)(42,22.5)
\bezier{60}(39,30)(40.5,30.75)(42,31.5)

\put(9,3){\line(-2,1){6.00}}
\put(9,3){\line(2,1){6.00}}
\put(21,3){\line(-2,1){6.00}}
\put(21,3){\line(2,1){6.00}}
\put(33,3){\line(-2,1){6.00}}
\put(33,3){\line(2,1){6.00}}

\put(3,6){\line(0,1){6.00}}
\put(15,6){\line(0,1){6.00}}
\put(27,6){\line(0,1){6.00}}
\put(39,6){\line(0,1){6.00}}

\put(9,15){\line(-2,-1){6.00}}
\put(9,15){\line(2,-1){6.00}}
\put(21,15){\line(-2,-1){6.00}}
\put(21,15){\line(2,-1){6.00}}
\put(33,15){\line(-2,-1){6.00}}
\put(33,15){\line(2,-1){6.00}}

\put(9,15){\line(0,1){6.00}}
\put(21,15){\line(0,1){6.00}}
\put(33,15){\line(0,1){6.00}}
\put(9,21){\line(-2,1){6.00}}
\put(9,21){\line(2,1){6.00}}
\put(21,21){\line(-2,1){6.00}}
\put(21,21){\line(2,1){6.00}}
\put(33,21){\line(-2,1){6.00}}
\put(33,21){\line(2,1){6.00}}
\put(3,24){\line(0,1){6.00}}
\put(15,24){\line(0,1){6.00}}
\put(27,24){\line(0,1){6.00}}
\put(39,24){\line(0,1){6.00}}
\put(9,33){\line(-2,-1){6.00}}
\put(9,33){\line(2,-1){6.00}}
\put(21,33){\line(-2,-1){6.00}}
\put(21,33){\line(2,-1){6.00}}
\put(33,33){\line(-2,-1){6.00}}
\put(33,33){\line(2,-1){6.00}}
\put(14,3.6){$\scriptstyle v_1$}
\put(9.3,1.0){$\scriptstyle v_2$}

\put(0,28.5){\line(2,1){15.00}}
\put(0,22.5){\line(2,1){27.00}}
\put(0,16.5){\line(2,1){39.00}}
\put(0,10.5){\line(2,1){39.00}}
\put(3,6){\line(2,1){39.00}}
\put(3,0){\line(2,1){39.00}}
\put(15,0){\line(2,1){27.00}}
\put(27,0){\line(2,1){15.00}}

\color{blue}
\put(14,1.2){$\scriptstyle \ga_1$}
\put(9.2,7.5){$\scriptstyle \ga_2$}
\put(9,3){\vector(1,0){12.0}}
\put(9,3){\vector(2,3){6}}
\end{picture}
\hspace{2mm}
\begin{picture}(25,30)(0,0)
\put(3,0){\emph{d})}
\put(10,3){\circle*{1.5}}
\put(22,19){\circle{1.5}}
\put(9,0){$v_2$}
\put(22,20){$v_1$}

\bezier{200}(10,3)(6,21)(22,19)
\bezier{200}(10,3)(27,2)(22,19)
\put(8,23){$\cG_*$}

\bezier{200}(10,3)(12,15)(22,19)
\bezier{200}(10,3)(21,8)(22,19)
\end{picture}

\begin{picture}(142,20)(0,0)
\put(0,3){\emph{c})}
\put(10,10){\circle{1.5}}
\put(20,5){\circle*{1.5}}
\put(30,10){\circle{1.5}}
\put(40,5){\circle*{1.5}}
\put(50,10){\circle{1.5}}
\put(60,5){\circle*{1.5}}
\put(70,10){\circle{1.5}}
\put(80,5){\circle*{1.5}}

\put(80,5){\line(2,1){5.00}}
\put(10,10){\line(-2,-1){5.00}}
\put(20,5){\line(2,1){10.00}}
\put(40,5){\line(2,1){10.00}}
\put(60,5){\line(2,1){10.00}}
\put(10,10){\line(2,-1){10.00}}
\put(30,10){\line(2,-1){10.00}}
\put(50,10){\line(2,-1){10.00}}
\put(70,10){\line(2,-1){10.00}}
\put(19,2){$v_2$}
\put(28,12){$v_1$}

\bezier{200}(10,10)(18,13)(20,5)
\bezier{200}(10,10)(12,2)(20,5)
\bezier{200}(30,10)(38,13)(40,5)
\bezier{200}(30,10)(32,2)(40,5)
\bezier{200}(50,10)(58,13)(60,5)
\bezier{200}(50,10)(52,2)(60,5)
\bezier{200}(70,10)(78,13)(80,5)
\bezier{200}(70,10)(72,2)(80,5)
\put(4,11){$\cG_\gt$}
\color{blue}
\put(20,5){\vector(1,0){20.00}}
\put(29,3){$\scriptstyle\ga_1$}
\end{picture}
\caption{\scriptsize\emph{a}) The diamond lattice $\cG$ is obtained by stacking together infinitely many copies of the hexagonal lattice $\bG$ along the height axis. The copies are connected in a periodic way by edges between white vertices in a lower copy and black vertices in the upper one. The blue vectors $\ga_1,\ga_2,\ga_3$ are the periods of $\cG$. \emph{b}) The subcovering $\cG_\gt$ of $\cG$ with the red chiral vector $\gt=\ga_1+\ga_2-\ga_3$. \emph{c}) The subcovering $\cG_\gt$ of $\cG$ with $\gt=\{\gt_1,\gt_2\}$, where the chiral vectors $\gt_1=\ga_1-\ga_3$ and $\gt_2=\ga_2-\ga_3$ (marked in green color). \emph{d}) The fundamental graph $\cG_*$ of the diamond lattice $\cG$ and any its subcovering $\cG_\gt$ with a primitive set $\gt$ of chiral vectors.} \label{Fcon}
\end{figure}

\begin{example} \lb{Exdl} The diamond lattice $\cG$ shown in Fig. \ref{Fcon}\emph{a} is a $\Z^3$-periodic graph, the vectors $\ga_1,\ga_2,\ga_3$ are the periods of $\cG$. The subcovering $\cG_\gt=\cG/\Z\gt$ of the diamond lattice $\cG$ with the chiral vector $\gt=(1,1,-1)$ is shown in Fig. \ref{Fcon}\emph{b}. Note that this subcovering of $\cG$ is isomorphic to the square lattice. Let $\G_\gt$ be the sublattice of $\Z^3$ with the basis $\gt:=\{\gt_1,\gt_2\}$, where $\gt_1=(1,0,-1)$ and $\gt_2=(0,1,-1)$. Then the subcovering $\cG_{\gt}=\cG/\G_\gt$ of the diamond lattice $\cG$ is the $\Z$-periodic graph with the period $\ga_1$ shown in Fig. \ref{Fcon}\emph{c}.
\end{example}

\begin{remark}\lb{Rade}
A specific class of subcovering graphs was studied in \cite{S24}. Namely, in \cite{S24} the author dealt with perturbations of a periodic graph by adding edges in a periodic way (without changing the vertex set). This perturbed periodic graph can be considered as a subcovering of a higher dimensional periodic graph obtained by stacking together infinitely many copies of the unperturbed one and connecting them by edges in an appropriate periodic way. For example, the square lattice perturbed by adding an edge $\{v,v+\ga_1+\ga_2\}$ at each vertex $v$ (i.e., the triangular lattice), see Fig. \ref{fig2}\emph{b}, is a subcovering of the cubic lattice (Fig.~\ref{fig2}\emph{a}). The cubic lattice is obtained by stacking together infinitely many copies of the square lattice and connecting them by edges $\{v,v+\ga_3\}$, $\forall\,v\in\Z^3$. The hexagonal lattice perturbed by adding edges between the white and black vertices as shown in Fig. \ref{Fcon}\emph{b} is a subcovering of the diamond lattice. The diamond lattice is obtained by stacking together infinitely many copies of the hexagonal lattice and connecting them by edges as shown in Fig. \ref{Fcon}\emph{a}.
\end{remark}

\subsection{Spectra of Schr\"odinger operators on subcovering graphs} \lb{Se1.3} In order to describe the spectrum of the Schr\"odinger operators on subcovering graphs we follow the approach which was used in \cite{KP07} for nanotubes.

Denote by $H_\gt=\D_\gt+Q$ the Schr\"odinger operator on the subcovering $\cG_\gt=\cG/\G_\gt$ with the potential $Q$ and the weight $\o$ induced by the $\Z^d$-periodic potential $Q$ and weight $\o$ on the original periodic graph $\cG$ (we use the same notation $Q$ and $\o$ for the induced potential and weight). One can think of $H_\gt$ as of the Schr\"odinger operator $H=\D+Q$ on $\cG=(\cV,\cE)$ acting on functions $f:\cV\to\C$ that are periodic with respect to the sublattice $\G_\gt\ss\Z^d$.

Recall that the Floquet-Bloch theory provides the direct integral expansion \er{die0} for the Schr\"odinger operator $H$ on the original periodic graph $\cG$. Since each function on $\cG_\gt$ lifts to a $\G_\gt$-periodic function $f$ on $\cG$, using \er{FlBc}, we obtain
\[\lb{qpco}
f(v+\gt_s)=e^{i\lan \gt_s,k\ran}f(v)=f(v),\qqq\forall\, v\in\cV, \qqq \forall\,s\in\N_{d_o},
\]
where $\gt_1,\ldots,\gt_{d_o}\in\Z^d$, $d_o<d$, are the chiral vectors of the subcovering $\cG_\gt$.

Let $\cT$ be the $d_o\ts d$ integer matrix whose rows are the chiral vectors $\gt_s$, $s\in\N_{d_o}$:
\[\lb{macT}
\cT=\left(
\begin{array}{ccc}
t_{11}& \ldots & t_{1d}\\
\ldots & \ldots & \ldots\\
t_{d_o1}& \ldots & t_{d_od}\\
\end{array}\right),\qqq \textrm{where}\qqq \gt_s=(t_{s1},\ldots,t_{sd})\in\Z^d.
\]
This matrix $\cT$ will be called the \emph{chiral matrix} of the subcovering $\cG_\gt$.

From \er{qpco} it follows that $\cT k\in2\pi\Z^{d_o}$. Then for the Schr\"odinger operator $H_\gt$ on the subcovering $\cG_\gt$ we have the direct integral decomposition
$$
H_\gt=\int^\oplus_{\cB_\gt}H(k)\,dk,
$$
where $H(k)$ is the Floquet operator for the original periodic graph $\cG$ defined by \er{Flop}, and $\cB_\gt$ is the following subset of the Brillouin zone $\cB=(-\pi,\pi]^d$ for $\cG$:
\[\lb{bgt}
\cB_\gt=\big\{k\in(-\pi,\pi]^d: \cT k\in2\pi\Z^{d_o}\big\}.
\]
This yields that the spectrum of $H_\gt$ is given by
\[\lb{spDt}
\s(H_\gt)=\bigcup_{k\in\cB_\gt}\s\big(H(k)\big),
\]
and the dispersion relation for $H_\gt$ is the restriction of the dispersion relation for $H$ to $\cB_\gt$.

From \er{spDt} it is clear that for any subcovering $\cG_\gt$ of a periodic graph $\cG$,
\[\lb{insp}
\s(H_\gt)\subseteq\s(H).
\]
In particular, if $\g$ is a gap of $H$, then the interval $\g$ is contained in a gap of $H_\gt$. Moreover, if $\L$ is a flat band of $H$, then $\L$ is also a flat band of $H_\gt$. But during the restriction of the dispersion relation for $H$ to $\cB_\gt$ new gaps might open and new flat bands might appear for $H_\gt$.

\medskip

The paper is organized as follows. In Section \ref{Sec2} we formulate our main results:

$\bullet$ we show that the spectra of the Schr\"odinger operators $H$ on a periodic graph $\cG$ and $H_\gt$ on any its subcovering $\cG_\gt$ with a \emph{primitive} set $\gt$ of chiral vectors (i.e., a set $\gt$ which can be completed to a basis of the lattice $\Z^d$) consist of the same number of spectral bands, and the corresponding band edges are asymptotically close to each other as the chiral vectors are long enough (Proposition \ref{PAig});

$\bu$ we formulate a simple criterion for the subcovering $\cG_\gt$ to be isospectral to the original periodic graph $\cG$ (Proposition \ref{PAig}.\emph{i});

$\bullet$ we obtain asymptotics of the band edges of the Schr\"odinger operator $H_\gt$ on the subcovering $\cG_\gt$ as the length of all chiral vectors from the primitive set $\gt$ tends to infinity (Theorem \ref{TAsy1}).

\medskip

\no We illustrate the obtained results by some examples of periodic graphs and their subcoverings. In particular, we describe all subcoverings of the $d$-dimensional lattice, the hexagonal lattice and the diamond lattice (with the unit edge weight) which are isospectral to the original periodic graphs (Examples \ref{Ehex}, \ref{ExSl} and \ref{Exa3}, respectively).

In Section \ref{Sec3} we prove our main results. The proofs are based on the connection between the dispersion relations for the Schr\"odinger operators $H$ on a periodic graph $\cG$ and $H_\gt$ on its subcovering $\cG_\gt$. In the proof we also use the crucial fact that the periodic graph $\cG$ and its subcovering $\cG_\gt$ with a primitive set $\gt$ of chiral vectors have the same fundamental graph. Section \ref{Sec4} is devoted to examples of periodic graphs and their asymptotically isospectral and just isospectral subcoverings.

\section{Main results}
\setcounter{equation}{0}
\lb{Sec2}
\subsection{Primitive sets in $\Z^d$} A set $\gt=\{\gt_1,\ldots,\gt_{d_o}\}\ss\Z^d$, where $d_o<d$, is said to be \emph{primitive} if $\gt$ is a basis for the lattice $\Z^d\cap\span(\gt)$ or, equivalently, if $\gt$ can be completed to a basis of $\Z^d$. Here $\span(\gt)$ is the set of all linear combinations of the vectors from $\gt$ with real coefficients. When $d_o=1$, the single vector $\gt_1$ in $\gt$ is called a \emph{primitive} vector. In this case we will often omit the subscript 1 in $\gt_1$ and just denote this single vector by $\gt$.

There is a useful characterization of primitive sets. Let $\cT$ be the $d_o\ts d$ integer matrix whose rows are the vectors $\gt_1,\ldots,\gt_{d_o}\in\Z^d$. It is known, see, e.g., \cite{C97}, that $\{\gt_1,\ldots,\gt_{d_o}\}\ss\Z^d$ is a primitive set if and only if the greatest common divisor of all the $d_o$-th order minors of the matrix $\cT$ is one. In particular, a vector $\gt=(t_1,\ldots,t_d)\in\Z^d$ is primitive if and only if its components $t_1,\ldots,t_d$ are coprime integers, i.e., $\gt$ is not an integral multiple of other vectors (other than $\pm\gt$) of $\Z^d$. This characterization of primitive sets is independent of the choice of the basis of $\Z^d$, i.e., if the rows of the matrix $\cT$ form a primitive set and $U$ is a $d\ts d$ unimodular matrix (which relates different bases of $\Z^d$), then the rows of the matrix $\cT U$ also form a primitive set.

\begin{remark}
\emph{i}) Primitive vectors of the lattice $\Z^d$ are also called \emph{visible}, since only they can be seen from the origin of the lattice. Non-primitive lattice points are hidden from view by other lattice points lying in the line of sight.

\emph{ii}) Each vector in a primitive set is a primitive vector. But not every set of primitive vectors is a primitive set.  For example, $\gt_1=(1,1,0)$ and $\gt_2=(1,-1,0)$ are primitive vectors of $\Z^3$, but $\{\gt_1,\gt_2\}$ is not a primitive set in $\Z^3$. Indeed, $\span(\gt_1,\gt_2)=\R^2\ts\{0\}$, so
$$
\Z^3\cap\span(\gt_1,\gt_2)=\Z^2\ts\{0\}.
$$
But $\{\gt_1,\gt_2\}$ is not a basis for the lattice $\Z^2\ts\{0\}$ (for example, the vector $(1,0,0)\in\Z^2\ts\{0\}$ does not belong to the lattice with the basis $\{\gt_1,\gt_2\}$).
\end{remark}

\subsection{Spectra of Schr\"odinger operators on subcoverings with primitive sets of chiral vectors}
Let $\cG$ be a $\Z^d$-periodic graph with the fundamental graph $\cG_*=\cG/\Z^d=(\cV_*,\cE_*)$. From now on we will consider only subcoverings $\cG_\gt$ of $\cG$ with primitive sets $\gt\ss\Z^d$ of chiral vectors $\gt_1,\ldots,\gt_{d_o}$, $d_o<d$.

We describe the spectrum of the Schr\"odinger operator $H_\gt$ on the subcovering $\cG_\gt$. Denote by $K_j^\pm$ the level sets corresponding to the band edges $\l_j^\pm$ of the Schr\"odinger operator $H$ on $\cG$, i.e.,
\[\lb{Kjpm}
K_j^\pm=\big\{k\in(-\pi,\pi]^d: \l_j(k)=\l_j^\pm\big\},\qqq j\in\N_\n, \qqq \n=\#\cV_*,
\]
where $\l_j(\cdot)$ are the band functions of $H$.

\begin{proposition}\lb{PAig}
Let $\cG$ be a $\Z^d$-periodic graph, and $\gt=\{\gt_1,\ldots,\gt_{d_o}\}\ss\Z^d$ be a primitive set. Then the spectrum of the Schr\"odinger operator $H_\gt=\D_\gt+Q$ on the subcovering $\cG_\gt$ has the form
\[\lb{debe}
\s(H_\gt)=
\bigcup_{j=1}^\n\s_j(H_\gt),\qqq
\s_j(H_\gt)=\big[\l_j^-(\gt),\l_j^+(\gt)\big]\subseteq\big[\l_j^-,\l_j^+\big],\qqq \n=\#\cV_*,
\]
where $\l_j^\pm$ are the band edges of the Schr\"odinger operator $H=\D+Q$ on $\cG$ defined in \er{ban.1H}. Moreover, the following statements hold true:

i) $\l_j^\pm(\gt)=\l_j^\pm$ if and only if
\[\lb{soeq}
\cT k_o\in2\pi\Z^{d_o}
\]
for some $k_o\in K_j^\pm$, where $K_j^\pm$ are given by \er{Kjpm}, and $\cT$ is the chiral matrix of $\cG_\gt$ defined by \er{macT}. In particular,

$\bu$ If $0\in K_j^\pm$, then $\l_j^\pm(\gt)=\l_j^\pm$.

$\bu$ If $(\pi,\ldots,\pi)\in K_j^\pm$ and the sum of the entries in each row of $\cT$ is even, then $\l_j^\pm(\gt)=\l_j^\pm$.

ii) Let $(\cG_\gt)$ be a sequence of the subcoverings of $\cG$ such that the smallest singular value $\t$ of their $d_o\ts d$ chiral matrices $\cT$ goes to infinity. Then
\[\lb{as00}
\lim\limits_{\t\to\iy}\l_j^\pm(\gt)=\l_j^\pm,\qqq \forall\,j\in\N_\n.
\]
\end{proposition}

\begin{remark}\lb{RAs}
\emph{i}) Proposition \ref{PAig} shows that the spectral bands
$$
\s_j(H)=\big[\l_j^-,\l_j^+\big]
\qqq\textrm{and}\qqq \s_j(H_\gt)=\big[\l_j^-(\gt),\l_j^+(\gt)\big],\qqq j\in\N_\n, \qqq \n=\#\cV_*,
$$
of the Schr\"odinger operators $H$ on a periodic graph $\cG$ and $H_\gt$ on its subcovering $\cG_\gt$ with a primitive set $\gt$ satisfy $\s_j(H_\gt)\subseteq\s_j(H)$.

\emph{ii}) From Proposition \ref{PAig} it also follows that a periodic graph $\cG$ is \emph{isospectral} to its subcovering $\cG_\gt$ with a primitive set $\gt$, if and only if for each $j\in\N_\n$ the condition \er{soeq} holds for some $k^-\in K_j^-$ and for some $k^+\in K_j^+$, where $K_j^\pm$ are given by \er{Kjpm}. By \emph{isospectrality} of periodic graphs we mean that the spectra of the Schr\"odinger operators on the graphs consist of the same number of bands and the corresponding bands coincide as sets.

\emph{iii}) If the set $\gt\ss\Z^d$ is not primitive, then Proposition \ref{PAig} is not true any more. For example, if $d_o=1$ and the single chiral vector $\gt=p\gt_o$ for some primitive $\gt_o\in\Z^d$ and integer $p\geq2$, then the fundamental graph of the subcovering $\cG_\gt$ has $p\n$ vertices, and, consequently, the spectrum of the Schr\"odinger operator $H_\gt$ on $\cG_\gt$ consists of $p\n$ spectral bands. Note that $\cG_{\gt_o}$ is a $p$-fold subcovering of $\cG_{p\gt_o}$, which yields that $\s(H_{\gt_o})\subseteq\s(H_{p\gt_o})$, see, e.g., \cite{A95}.

\emph{iv}) Since $\gt_1,\ldots,\gt_{d_o}$ are linearly independent, the matrix $\cT\cT^\top$, where $\cT^\top$ is the transpose of $\cT$, is positive definite. The smallest singular value $\t$ of $\cT$ is the smallest eigenvalue of $(\cT\cT^\top)^{1/2}$. The condition $\t\to\iy$ implies that the Euclidean norm $|\gt_s|$ of each chiral vector $\gt_s\in\Z^d$ tends to infinity. Indeed, the diagonal entries of the matrix $\cT\cT^\top$ are $|\gt_1|^2,\ldots,|\gt_{d_o}|^2$. Then the smallest eigenvalue $\t^2$ of $\cT\cT^\top$ satisfies
$$
\t^2\leq\min\big\{|\gt_1|^2,\ldots,|\gt_{d_o}|^2\big\},
$$
see, e.g., Theorem 4.3.15 in \cite{HJ85}. This yields that
\[\lb{gtin}
|\gt_s|\to\iy,\qqq s\in\N_{d_o}, \qqq\textrm{as}\qqq  \t\to\iy.
\]
Note that if $d_o=1$, then $\t=|\gt_1|$. The condition \er{gtin} in turn implies that when $\t\to\iy$

$\bu$  the "radius" of the cylinder $\R^d/\G_\gt$, in which the subcovering $\cG_\gt$ is naturally embedded (see Section \ref{Sscg}), tends to infinity (recall that $\G_\gt$ is a sublattice of $\Z^d$ with the basis $\gt=\{\gt_1,\ldots,\gt_{d_o}\}$);

$\bu$ the set $\cB_\gt$ defined by \er{bgt} is dense in the Brillouin zone $\cB=(-\pi,\pi]^d$.

\emph{v}) Due to \er{as00}, the band edges $\l_j^\pm(\gt)$ of $H_\gt$ are asymptotically close to the corresponding band edges $\l_j^\pm$ of $H$ as $\t\to\iy$. Then we say that $\cG_\gt$ is \emph{asymptotically isospectral} to $\cG$ as $\t\to\iy$. Notion of asymptotic isospectrality for Schr\"odinger operators on periodic graphs was introduced in \cite{S24} when studying the behavior of the spectrum of the Schr\"odinger operators on periodic graphs under perturbations of graphs by adding edges (between existing vertices) in a periodic way.

\emph{vi}) Since the spectrum of the Schr\"odinger operators on periodic graphs has band structure, there exist some notions of isospectrality for periodic graphs. The problem of Floquet and Fermi isospectrality for the discrete Schr\"odinger operators on the lattice $\Z^d$ was studied in \cite{GKT93,K89,L23} and \cite{L24}, respectively.
\end{remark}

\subsection{Asymptotics of the band edges.} In this subsection we present asymptotics of the band edges $\l_j^\pm(\gt)$ of the Schr\"odinger operator $H_\gt$ on the subcovering $\cG_\gt$ as all chiral vectors $\gt_1,\ldots,\gt_{d_o}$ from the set $\gt$ are long enough.

From now on we impose the following assumptions on a band function $\l_j(\cdot)$ of the Schr\"odinger operator $H$ on the periodic graph $\cG$:

\medskip

\textbf{Assumption A}

\textbf{A1} The band function $\l_j(\cdot)$ of $H$ has a non-degenerate minimum $\l_j^-$ (maximum $\l_j^+$) at some point $k_o\in\cB=(-\pi,\pi]^d$, and this extremum $\l_j^\pm$ is attained by the single band function $\l_j(\cdot)$.

\textbf{A2} $k_o$ is the only (up to evenness) minimum (maximum) point of $\l_j(\cdot)$ in $\cB$, see Remark \ref{prbf}.\emph{ii}.

\begin{theorem}\lb{TAsy1}
Let $\gt=\{\gt_1,\ldots,\gt_{d_o}\}\ss\Z^d$ be a primitive set, and let for some $j\in\N_\n$ the band function $\l_j(\cdot)$ of the Schr\"odinger operator $H$ on a $\Z^d$-periodic graph $\cG$ satisfy the Assumption A. Then the lower band edge $\l_j^-(\gt)$ (the upper band edge $\l_j^+(\gt)$) of the Schr\"odinger operator $H_\gt$ on the subcovering $\cG_\gt$ of $\cG$ has the following asymptotics
\[\lb{aslk1}
\begin{array}{l} \l_j^\pm(\gt)=\l_j^\pm\mp\frac12\,\big|(\cT\bH^{-1}\cT^\top)^{-1/2}\bx_o\big|^2+g(\gt),
\qqq \textrm{where}\\[8pt]
\bH=\mp\mathrm{Hess}\,\l_j(k_o),\qqq \bx_o=(\cT k_o\hspace{-3mm}\mod2\pi\Z^{d_o})\in(-\pi,\pi]^{d_o},\\[8pt]
g(\gt)=O\big(\frac1{\t^3}\big),\qqq \textrm{as}\qqq \t\to\iy.
\end{array}
\]
Here $\mathrm{Hess}\,\l_j(k_o)$ is the Hessian of $\l_j$ at $k_o$; $\cT$ is the chiral matrix of $\cG_\gt$ defined by \er{macT}; $\t$ is the smallest singular value of $\cT$, and $|\cdot|$ denotes the Euclidean norm of a vector.

In particular, if $\cG_\gt$ is a subcovering of $\cG$ with a single primitive vector $\gt\in\Z^d$, then
\[\lb{aslk}
\begin{array}{l}
\displaystyle\l_j^\pm(\gt)=\l_j^\pm\mp\frac12\cdot\frac{|x_o|^2}{|\bH^{-1/2}\,\gt|^2}
+g(\gt),\qqq\textrm{where}\\[12pt]
x_o=(\lan\gt,k_o\ran\hspace{-3mm}\mod2\pi)\in(-\pi,\pi],\\[8pt]
g(\gt)=\textstyle O\big(\frac1{\t^3}\big),\qqq \textrm{as}\qqq \t=|\gt|\to\iy.
\end{array}
\]
Moreover, if all components of $k_o$ are equal to either 0 or $\pi$, then $g(\gt)=O\big(\frac1{\t^4}\big)$ in \er{aslk1} and \er{aslk}.
\end{theorem}

\begin{remark}\lb{RAsi}
\emph{i}) Assumption A1 is often used to establish many important properties of periodic operators such as electron's effective masses in solid state physics \cite{AM76}, Green's function asymptotics \cite{KKR17} and other issues. It is an old and widely believed conjecture in mathematical physics that generically (with respect to perturbations of the periodic potential and edge weights) Assumption A1 holds true. For recent progress in proving the band edges nondegeneracy conjecture see \cite{FS24} and references therein.

\emph{ii}) The asymptotics \er{aslk1} can be carried over to the case when the band edge $\l_j^\pm$ occurs at finitely many quasimomenta $k_o\in\cB$ (instead of assuming the condition A2) by taking minimum (maximum) among the asymptotics coming from all these non-degenerate isolated extrema.

\emph{iii}) The second term in the asymptotics \er{aslk1} is equal to zero if and only if $\cT k_o\in2\pi\Z^{d_o}$. This agrees with Proposition~\ref{PAig}, see the condition \er{soeq}.

\emph{iv}) For $\t$ large enough, the error term $g(\gt)$ in \er{aslk1} and \er{aslk} satisfies $|g(\gt)|\leq C\t^{-3}$, where $C$ is a positive constant depending on the Hessian of the band function $\l_j$ at the extrema point $k_o$ and the number $d_o$ of the chiral vectors but not depending on the chiral matrix $\cT$ (i.e., not depending on the choice of the subcovering $\cG_\gt$).

\emph{v}) For a particular case when $d_o=1$ and the single chiral vector has the specific form $\gt=(t_1,\ldots,t_{d-1},-1)\in\Z^d$, i.e., $t_d=-1$, the statements of Proposition~\ref{PAig} and Theorem \ref{TAsy1} were obtained in \cite{S24}, see also Remark \ref{Rade}.
\end{remark}

\subsection{Examples.} In this section we illustrate the obtained results by some simple examples of periodic graphs and their subcoverings. In all examples we will consider the Schr\"odinger operator $H$ defined by \er{Sh}, \er{ALO}
with the edge weight $\o\equiv1$. The proofs of the examples are based on Proposition~\ref{PAig} and Theorem \ref{TAsy1} and given in Section \ref{Sec4}.

\medskip

First we consider the hexagonal lattice $\bG$ which is a $\Z^2$-periodic graph with the periods $\ga_1,\ga_2$, see  Fig.~\ref{fig1}\emph{a}. The fundamental graph $\bG_*$ of $\bG$ consists of two vertices $v_1$ and $v_2$ and three multiple edges connecting these vertices, see Fig.~\ref{fig1}\emph{d}. Let $H=\D+Q$ be the Schr\"odinger operator on $\bG$ with a $\Z^2$-periodic potential $Q$. Without loss of generality (add a constant to $Q$ if necessary) we may assume that the potential $Q$ is given by
\[\lb{poGL}
Q(v_1)=-Q(v_2)=:q>0.
\]
Then the spectrum of $H$ has the form
$$
\s(H)=\big[3-\sqrt{9+q^2}\,,3-q\big]\cup\big[3+q,3+\sqrt{9+q^2}\,\big],
$$
see, e.g., Lemma 4.1 in \cite{S24}. We describe the spectrum of the Schr\"odinger operator $H_\gt=\D_\gt+Q$ on subcoverings (nanotubes) $\bG_\gt$ of the hexagonal lattice $\bG$.

\begin{example}\lb{Ehex}
Let $\bG_\gt$ be a nanotube with a primitive chiral vector $\gt=(t_1,t_2)\in\Z^2$. Then the spectrum $\s(H_\gt)$ of the Schr\"odinger operator $H_\gt=\D_\gt+Q$ on $\bG_\gt$ with the periodic potential $Q$ defined by \er{poGL} has the form
\[\lb{stin}
\s(H_\gt)=\big[3-\sqrt{9+q^2}\,,3-q-\ve(\gt)\big]\cup
\big[3+q+\ve(\gt),3+\sqrt{9+q^2}\,\big].
\]

\emph{i}) If $t_1-t_2\in3\Z$, then $\ve(\gt)=0$.

\emph{ii}) If $t_1-t_2\notin3\Z$, then
$$
\ve(\gt)=\frac{\pi^2}{6q(t_1^2+t_1t_2+t_2^2)}+\textstyle O\big(\frac1{|\gt|^3}\big),\qqq \textrm{as} \qqq |\gt|\to\infty.
$$
\end{example}

\begin{remark} \emph{i}) For any primitive $\gt\in\Z^2$, the nanotube $\bG_\gt$ is asymptotically isospectral to the hexagonal lattice $\bG$ as $|\gt|\to\iy$. If $t_1-t_2\in3\Z$, then the Schr\"odinger operators $H$ on $\bG$ and $H_\gt$ on $\bG_\gt$ have the same spectrum, i.e., $\bG$ and $\bG_\gt$ are just isospectral. If $\gt$ is not primitive, i.e., $\gt=p\gt_o$ for some primitive $\gt_o\in\Z^2$ and integer $p\geq2$, then the spectrum of the Schr\"odinger operator $H_\gt$ on $\bG_\gt$ consists of $2p$ spectral bands, some of which may be degenerate. But, due to \er{insp}, the inclusion $\s(H_{p\gt_o})\subseteq\s(H)$ still holds.

For example, if $\bG_\gt$ is the zig-zag nanotube with the non-primitive chiral vector $\gt=(2,0)$ (see Fig.~\ref{fig1}\emph{c}), then the spectrum of the Schr\"odinger operator $H_\gt=\D_\gt+Q$ on $\bG_\gt$ with the periodic potential $Q$ defined by \er{poGL} has the form
$$
\begin{array}{l}
\s(H_\gt)=\s_{ac}(H_\gt)\cup\s_{fb}(H_\gt),\qqq \s_{fb}(H_\gt)=\{3\pm\sqrt{1+q^2}\,\},\\ [8pt] \s_{ac}(H_\gt)=\big[3-\sqrt{9+q^2}\,,3-\sqrt{1+q^2}\,\big]\cup
\big[3+\sqrt{1+q^2}\,,3+\sqrt{9+q^2}\,\big],
\end{array}
$$
see, e.g., \cite{KK10}. The spectrum $\s(H_\gt)$ consists of four spectral bands (since the fundamental graph of $\bG_{\gt}$ has four vertices, see Fig.~\ref{fig1}\emph{e}), two of which are flat bands.

\emph{ii}) For any primitive $\gt\in\Z^2$, the spectrum of the Schr\"odinger operator $H_\gt$ on the nanotube $\bG_\gt$ is symmetric with respect to the point 3 and has a gap of length $2(q+\ve(\gt))>0$.

\emph{iii}) The spectra of discrete (and differential) Schr\"odinger operators on nanotubes with any chiral vector $\gt$ (including non-primitive ones) were studied in \cite{KP07}. Zig-zag ($\gt=(n,0)$) and armchair ($\gt=(n,n)$) nanotubes in magnetic and electric fields were considered in \cite{KK10}.
\end{remark}

As the second example we consider the $d$-dimensional lattice $\dL^d$, for $d=3$ see Fig.~\ref{fig2}\emph{a}. The lattice $\dL^d$ is a $\Z^d$-periodic graph with the periods $\ga_1,\ldots,\ga_d$. The spectrum of the Laplacian $\D$ on $\dL^d$ has the form $\s(\D)=[0,4d]$. We describe the Laplacian spectrum on the subcovering $\dL^d_\gt$ of the lattice $\dL^d$ with a primitive set $\gt$ of chiral vectors $\gt_1,\ldots,\gt_{d_o}\in\Z^d$, $d_o<d$  (for $d=3$, $d_o=1$ and $\gt_1=(1,1,-1)$ see Fig.~\ref{fig2}\emph{b}; for $d=3$, $d_o=2$ and $\gt_1=(1,1,-1)$, $\gt_2=(2,1,0)$ see Fig.~\ref{fig2}\emph{c}).

\begin{example}\lb{ExSl}
Let $\gt=\{\gt_1,\ldots,\gt_{d_o}\}\ss\Z^d$ be a primitive set of vectors
$$
\gt_s=(t_{s1},\ldots,t_{sd})\in\Z^d,\qqq s\in\N_{d_o}, \qqq d_o<d,
$$
and let $\dL^d_\gt$ be the corresponding subcovering of the $d$-dimensional lattice $\dL^d$. Then the spectrum $\s(\D_\gt)$ of the Laplacian $\D_\gt$ on $\dL^d_\gt$ is given by
\[\lb{sDLt}
\s(\D_\gt)=[0,\l_1^+(\gt)].
\]

\emph{i}) If $\varrho_s:=t_{s1}+\ldots+t_{sd}$ is even for all $s\in\N_{d_o}$, then $\l_1^+(\gt)=4d$.

\emph{ii}) If $\vr_s$ is odd for some $s\in\N_{d_o}$,  then
\[\lb{ase0}
\begin{array}{l}
\l_1^+(\gt)=4d-\pi^2\big|
(\cT\cT^\top)^{-1/2}\vr\big|^2+O\big(\frac1{\t^4}\big),\qq \textrm{as}\qq \t\to\iy,\\[6pt]
\textrm{where}\qqq 0\neq\vr=\big((\vr_1,\ldots,\vr_{d_o})\hspace{-3mm}\mod2\Z^{d_o}\big)\in\{0,1\}^{d_o},
\end{array}
\]
$\cT$ is the chiral matrix of $\dL^d_\gt$ defined by \er{macT}, and $\t$ is the smallest singular value of $\cT$.

In particular, if $\dL^d_\gt$ is a subcovering of $\dL^d$ with a single primitive vector $\gt=(t_1,\ldots,t_d)\in\Z^d$ and $t_1+\ldots+t_d$ is odd, then
\[\lb{ase1}
\textstyle\l_1^+(\gt)=4d-\dfrac{\pi^2}{\t^2}+O\big(\frac1{\t^4}\big),
\qqq \textrm{as}\qqq \t=|\gt|\to\iy.
\]
\end{example}

\begin{remark} \emph{i}) For any primitive set $\gt\ss\Z^d$, the $d$-dimensional lattice $\dL^d$ and its subcovering $\dL^d_\gt$ are asymptotically isospectral as $\t\to\iy$. If $\vr_s$ is even for all $s\in\N_{d_o}$, then $\l_1^+(\gt)$ does not depend on $\gt$ and the Laplacians on $\dL^d$ and $\dL^d_\gt$ have the same spectrum $[0,4d]$, i.e., $\dL^d$ and $\dL^d_\gt$ are just isospectral.

\emph{ii}) The triangular lattice shown in Fig.~\ref{fig2}\emph{b} is the subcovering $\dL^3_{\gt}$ of the cubic lattice $\dL^3$ with the chiral vector $\gt=(1,1,-1)$. It is known that the spectrum of the Laplacian $\D_{\gt}$ on the triangular lattice has the form $\s(\D_{\gt})=[0,9]$. On the other hand, since $\gt$ is primitive and the sum of its components is odd, then, due to Example \ref{ExSl}, we have $\s(\D_{\gt})=[0,\l_1^+(\gt)]$, and the asymptotics \er{ase1} (which is very rude in this case) yields
$$
\textstyle \l_1^+(\gt)\approx4\cdot3-\frac{\pi^2}{|\gt|^2}\approx8{.}71.
$$

\emph{iii}) To demonstrate the asymptotics \er{ase0} we consider the subcovering $\dL^3_\gt$ of the cubic lattice $\dL^3$ with the primitive set $\gt$ of the chiral vectors $\gt_1=(1,5,-1)$ and $\gt_2=(4,1,0)$. Since the sum of the components of $\gt_1$ ($\gt_2$) is odd, then, due to Example \ref{ExSl}, we have $\s(\D_\gt)=[0,\l_1^+(\gt)]$, and
$$
\l_1^+(\gt)\approx4\cdot3-\pi^2\big|
(\cT\cT^\top)^{-1/2}\vr\big|^2, \qq \textrm{where} \qq
\cT=\left(
\begin{array}{ccr}
1 & 5 & -1 \\
4 & 1 & 0
\end{array}\right),\qq
\vr=\left(
\begin{array}{c}
1 \\ 1
\end{array}\right).
$$
After some simple calculations, we obtain
$$
\textstyle\l_1^+(\gt)\approx12-\frac{13}{189}\,\pi^2\approx11{.}32,
$$
that agrees well with the spectrum $\s(\D_\gt)=[0,11.34]$ calculated numerically. Note that the smallest singular value of the matrix $\cT$ is $\t\approx3{.}42$.
\end{remark}

Finally, we consider the diamond lattice $\cG$ which is a $\Z^3$-periodic graph with the periods $\ga_1,\ga_2,\ga_3$, see Fig. \ref{Fcon}\emph{a}. The fundamental graph $\cG_*$ of $\cG$ consists of two vertices $v_1$ and $v_2$ and four multiple edges connecting these vertices, see Fig.~\ref{Fcon}\emph{d}. Let $H=\D+Q$ be the Schr\"odinger operator on $\cG$ with a $\Z^3$-periodic potential $Q$. Without loss of generality we may assume that the potential $Q$ is defined by
\[\lb{poDL}
Q(v_1)=-Q(v_2)=:q>0.
\]
Then the spectrum of $H$ has the form
\[\lb{spDL}
\s(H)=\big[4-\sqrt{16+q^2}\,,4-q\big]\cup\big[4+q,4+\sqrt{16+q^2}\,\big],
\]
see, e.g., Lemma 4.2 in \cite{S24}. We describe all subcoverings $\cG_\gt$ of the diamond lattice $\cG$ isospectral to $\cG$ (as well as to each other).

\begin{example}\lb{Exa3}
Let $\cG_\gt$ be a subcovering of the diamond lattice $\cG$ with a primitive set $\gt\ss\Z^3$ of chiral vectors, and let $H_\gt=\D_\gt+Q$ be the Schr\"odinger operator on $\cG_\gt$ with the periodic potential $Q$ defined by \er{poDL}. Then $\s(H_\gt)=\s(H)$ if and only if

$\bu$ $\gt$ consists of a single primitive vector, \emph{or}

$\bu$ $\gt=\{\gt_1,\gt_2\}$, where $\gt_s=(t_{s1},t_{s2},t_{s3})\in\Z^3$, $s=1,2$, such that for some $i,j=1,2,3$, $i\neq j$, the sum $t_{si}+t_{sj}$ is even for each $s=1,2$.\\
Here $\s(H)$ is given by \er{spDL}.
\end{example}

\begin{remark}
\emph{i}) The subcovering $\cG_\gt$ of the diamond lattice $\cG$ with the primitive chiral vector $\gt=(1,1,-1)$ shown in Fig.~\ref{Fcon}\emph{b} is isospectral to $\cG$.

\emph{ii}) If for a primitive set $\gt=\{\gt_1,\gt_2\}$ the condition of Example \ref{Exa3} does not hold, i.e., there are no $i,j=1,2,3$, $i\neq j$, such that the sum $t_{si}+t_{sj}$ is even for each $s=1,2$, then the spectrum of the Schr\"odinger operator $H_\gt$ on the subcovering $\cG_\gt$ of the diamond lattice has the form
\[\lb{dlni}
\s(H_\gt)=\big[4-\sqrt{16+q^2}\,,4-q-\ve(\gt)\,\big]
\cup\big[4+q+\ve(\gt)\,,4+\sqrt{16+q^2}\,\big],
\]
where $\ve(\gt)>0$ and $\ve(\gt)\to 0$ as $\t\to\iy$. Recall that $\t$ is the smallest singular value of the chiral matrix $\cT$ of $\cG_\gt$ defined by \er{macT}. For example, if $\gt_1=(1,0,-1)$ and $\gt_2=(0,1,-1)$, then the spectrum of $H_\gt$ on the subcovering $\cG_\gt$ with the primitive set $\gt=\{\gt_1,\gt_2\}$ shown in Fig.~\ref{Fcon}\emph{c} is given by \er{dlni}, where $\ve(\gt)=\sqrt{4+q^2}-q$.
\end{remark}

\section{Proof of the main results}\lb{Sec3}
\setcounter{equation}{0}

\subsection{A particular choice of the Brillouin zone.} Let $\cG$ be a $\Z^d$-periodic graph with the periods $\ga_1,\ldots,\ga_d$. Recall that the band functions $\l_j(k)$ of the Schr\"odinger operator $H$ on $\cG$ are $2\pi\Z^d$-periodic in $k\in\R^d$. In order to prove our main results we need a particular choice of the Brillouin zone for $\cG$ instead of $\cB=(-\pi,\pi]^d$. By a Brillouin zone we understand any choice of the fundamental domain of the lattice $2\pi\Z^d$ in the momentum space.

Let $\cG_\gt$ be a subcovering of $\cG$ with a primitive set $\gt=\{\gt_1,\ldots,\gt_{d_o}\}\ss\Z^d$, $d_o<d$. This primitive set $\gt$ can be extended to a basis
\[\lb{neba}
\{\gt_1,\ldots,\gt_{d_o},\gt_{d_o+1},\ldots,\gt_d\}
\]
of $\Z^d$. We denote by $\wt\cT$ the $d\ts d$ \emph{unimodular} matrix (i.e., the integer matrix with the determinant $\pm1$) whose rows are the vectors $\gt_1,\ldots\gt_d$:
\[\lb{macT1}
\wt\cT=(t_{ij})_{i,j=1}^d, \qq\textrm{where}\qq \gt_i=(t_{i1},\ldots,t_{id})\in\Z^d, \qq i\in\N_d,\qq \det\wt\cT=\pm1.
\]
Note that the first $d_o$ rows of $\wt\cT$ form the chiral matrix $\cT$ of the subcovering $\cG_\gt$, see \er{macT}. The change of the basis of $\Z^d$ (from $\{\ga_1,\ldots,\ga_d\}$ to the basis \er{neba}) determines the linear change in the quasimomentum~$k$:
$$
\k=\wt\cT k,\qqq k\in\R^d,
$$
and the new Brillouin zone $\wt\cB$ for $\cG$ given by
\[\lb{nBz}
\wt\cB=\big\{k\in\R^d :\wt\cT k\in(-\pi,\pi]^d\big\}.
\]

To describe the spectrum of the Schr\"odinger operator $H_\gt$ on the subcovering $\cG_\gt$ of the periodic graph $\cG$ we need to define the ranges of the band functions $\l_j(k)$ of the Schr\"odinger operator $H$ on $\cG$ restricted to the set $\cB_\gt$ given by
\[\lb{bgt1}
\cB_\gt=\big\{k\in(-\pi,\pi]^d: \cT k\in2\pi\Z^{d_o}\big\},
\]
where $\cT$ is the chiral matrix of the subcovering $\cG_\gt$, see Section \ref{Se1.3}.

\begin{lemma}\lb{Lrbf}
The ranges of the band functions $\l_j(k)$ of the Schr\"odinger operator $H$ on the periodic graph $\cG$, restricted to the set $\cB_\gt$ defined by \er{bgt1}, have the form
\[\lb{wtcb}
\begin{array}{l}
\displaystyle \l_j(\cB_\gt)=\l_j(\wt\cB_\gt)=
\big[\l_j^-(\gt),\l_j^+(\gt)\big],\qqq j\in\N_\n, \qqq \textrm{where}\\[6pt]
\displaystyle\l_j^-(\gt)=\min_{k\in\wt\cB_\gt}\l_j(k)=
\min_{k\in\cB_\gt}\l_j(k),\qqq
\l_j^+(\gt)=\max_{k\in\wt\cB_\gt}\l_j(k)=\max_{k\in\cB_\gt}\l_j(k),
\end{array}
\]
and the subset $\wt\cB_\gt$ of the new Brillouin zone $\wt\cB$, see \er{nBz}, is given by
\[\lb{nBzt}
\wt\cB_\gt=\big\{k\in\wt\cB :\cT k=0\big\}.
\]
\end{lemma}

\no \textbf{Proof.} The band functions $\l_j(k)$ are $2\pi\Z^d$-periodic in $k\in\R^d$. Taking into account that the set $\cS_\gt:=\{k\in\R^d: \cT k\in2\pi\Z^{d_o}\}$ is $2\pi\Z^d$-invariant, we obtain
$$
\l_j(\cB_\gt)=\l_j(\cS_\gt)=\l_j(\cS_\gt\cap\wt\cB)=\l_j(\wt\cB_\gt).
$$
The set $\wt\cB_\gt$ is connected, since its image under the linear isomorphism $\wt\cT$ defined by \er{macT1} is the connected set $\{0\}^{d_o}\ts(-\pi,\pi]^{d-d_o}$. Then, using that $\l_j(k)$ are continuous, we get
$$
\l_j(\wt\cB_\gt)=
\big[\l_j^-(\gt),\l_j^+(\gt)\big],\qqq j\in\N_\n,
$$
where $\l_j^\pm(\gt)$ are given in \er{wtcb}.\qq $\Box$

\medskip

The following example illustrates Lemma \ref{Lrbf}.

\begin{figure}[t!]\centering
\unitlength 0.8mm 
\linethickness{0.4pt}
\begin{picture}(150,95)

\put(0,45){\vector(1,0){150.00}}
\put(75,0){\vector(0,1){95.00}}
\multiput(60,30)(0,4){8}{\line(0,1){2}}
\multiput(60,30)(4,0){8}{\line(1,0){2}}
\multiput(90,30)(0,4){8}{\line(0,1){2}}
\multiput(60,60)(4,0){8}{\line(1,0){2}}
\put(142.5,0){\vector(-3,2){142.5}}
\put(0,82.5){\vector(2,-1){150}}

\bezier{180}(90,30)(45,60)(0,90)
\bezier{180}(92,29)(47,59)(2,89)
\bezier{180}(94,28)(49,58)(4,88)
\bezier{180}(96,27)(51,57)(6,87)
\bezier{180}(98,26)(53,56)(8,86)
\bezier{180}(100,25)(55,55)(10,85)
\bezier{180}(102,24)(57,54)(12,84)
\bezier{180}(104,23)(59,53)(14,83)
\bezier{180}(106,22)(61,52)(16,82)
\bezier{180}(108,21)(63,51)(18,81)
\bezier{180}(110,20)(65,50)(20,80)
\bezier{180}(112,19)(67,49)(22,79)
\bezier{180}(114,18)(69,48)(24,78)
\bezier{180}(116,17)(71,47)(26,77)
\bezier{180}(118,16)(73,46)(28,76)
\bezier{180}(120,15)(75,45)(30,75)
\bezier{180}(122,14)(77,44)(32,74)
\bezier{180}(124,13)(79,43)(34,73)
\bezier{180}(126,12)(81,42)(36,72)
\bezier{180}(128,11)(83,41)(38,71)
\bezier{180}(130,10)(85,40)(40,70)
\bezier{180}(132,9)(87,39)(42,69)
\bezier{180}(134,8)(89,38)(44,68)
\bezier{180}(136,7)(91,37)(46,67)
\bezier{180}(138,6)(93,36)(48,66)
\bezier{180}(140,5)(95,35)(50,65)
\bezier{180}(142,4)(97,34)(52,64)
\bezier{180}(144,3)(99,33)(54,63)
\bezier{180}(146,2)(101,32)(56,62)
\bezier{180}(148,1)(103,31)(58,61)
\bezier{180}(150,0)(105,30)(60,60)

\put(82,52.5){$\cB$}
\put(53,56.5){$\wt\cB$}
\put(91,46){\small$\pi$}
\put(53,42){\small$-\pi$}
\put(72,61){\small$\pi$}
\put(75,27){\small $-\pi$}
\put(147,40){$k_1$}
\put(69,92){$k_2$}
\put(146,10){$\k_1$}
\put(-4,92){$\k_2$}

\color{red}
\put(0,65){\line(3,-2){97.5}}
\put(0,75){\line(3,-2){112.5}}
\put(0,85){\line(3,-2){127.5}}
\put(0.3,95){\line(3,-2){142.5}}
\put(15,95){\line(3,-2){135.00}}
\put(30,95){\line(3,-2){120.00}}
\put(45,95){\line(3,-2){105.00}}
\put(-9,65){$L_{-3}$}
\put(-9,75){$L_{-2}$}
\put(-9,85){$L_{-1}$}
\put(5,92){$L_0$}
\put(20,92){$L_1$}
\put(35,92){$L_2$}
\put(50,92){$L_3$}
\put(100,68){$\ldots$}
\put(35,23){$\ldots$}
\end{picture}

\caption{\scriptsize Two choices of the Brillouin zone in terms of the quasimomenta $k=(k_1,k_2)\in(-\pi,\pi]^2$ (the square $\cB$ bounded by the dashed lines) and $\k=\wt\cT k$, $\k=(\k_1,\k_2)\in(-\pi,\pi]^2$ (the shaded parallelogram $\wt\cB$). The lines $L_n$ (marked in red color) are given by the equations $t_1k_1+t_2k_2=2\pi n$, $n\in\Z$, when $(t_1,t_2)=(2,3)$. The lines $L_0,L_{\pm1},L_{\pm2}$ intersect the square Brillouin zone $\cB$. The line $L_0$ (which coincides with the axis $\k_2$) is the single line among $L_n$ that enters the shaded  Brillouin zone $\wt\cB$.}
\label{fig0}
\end{figure}

\begin{example}
We consider the hexagonal lattice $\bG$, which is a $\Z^2$-periodic graph with the periods $\ga_1,\ga_2$, see Fig.~\ref{fig1}\emph{a}. Let $\bG_{\gt_1}$ be the subcovering of $\bG$ (the nanotube) with the primitive chiral vector $\gt_1=(2,3)$. The primitive vector $\gt_1=(2,3)$ can be completed to the basis $\{\gt_1,\gt_2\}$ of $\Z^2$ by the vector $\gt_2=(1,2)$. Thus, the matrix $\wt\cT$ defined by \er{macT1} has the form
$$
\wt\cT=\left(\begin{array}{cc}
2 & 3\\
1 & 2 \end{array}\right),\qqq \det\wt\cT=1.
$$
In Fig. \ref{fig0} we show two choices of the Brillouin zone for the hexagonal lattice $\bG$:\\
$\bu$ the square $\cB$ (bounded by the dashed lines) drawn in terms of the quasimomentum $k=(k_1,k_2)\in(-\pi,\pi]^2$, and \\
$\bu$ the parallelogram $\wt\cB$ (shaded in the figure) drawn in terms of the quasimomentum $\k=\wt\cT k$, $\k=(\k_1,\k_2)\in(-\pi,\pi]^2$.\\
The lines $L_n$ shown in the figure are given by the equations $t_1k_1+t_2k_2=2\pi n$, $n\in\Z$, where $\gt_1=(t_1,t_2)=(2,3)$ is the chiral vector of the nanotube $\bG_{\gt_1}$. The set $\cB_{\gt_1}$  defined by \er{bgt1} consists of the five segments of the lines $L_0,L_{\pm1},L_{\pm2}$ inside the square Brillouin zone $\cB$. The set $\wt\cB_{\gt_1}$ defined by \er{nBzt} consists of only the part of the line $L_0$ inside the Brillouin zone $\wt\cB$. The line $L_0$ coincides with the axis $\k_2$.

The dispersion relation for the Schr\"odinger operator $H=\D+Q$ on $\bG$ with the unit edge weight and a $\Z^2$-periodic potential $Q$ defined by \er{poGL} is plotted in Fig.~\ref{fDR} (left and center), when $q=1$. The Brillouin zone is parameterized by the quasimomentum $k$ in the left figure and by the quasimomentum $\k=\wt\cT k$ in the center figure. The band functions of the Schr\"odinger operator $H_{\gt_1}$ on the nanotube $\bG_{\gt_1}$ are the restrictions of the dispersion relation for $H$ to the segment of the line $L_0$ given by the equation $2k_1+3k_2=0$ (or, $\k_1=0$, in terms of $\k$) inside the Brillouin zone $\wt\cB$, see Fig.~\ref{fDR} (right). The ranges of the band functions of $H$ and $H_{\gt_1}$ (i.e., the spectra of $H$ and $H_{\gt_1}$) have the form
$$
\s(H)=[3-\sqrt{10}\,,2]\cup[4,3+\sqrt{10}\,],\qqq
\s(H_{\gt_1})\approx[3-\sqrt{10}\,,1{.}92]\cup[4{.}08,3+\sqrt{10}\,].
$$
They are close to each other but do not coincide.
\end{example}

\begin{figure}[t!]\centering
\includegraphics[width=4.5cm,height=4.5cm]{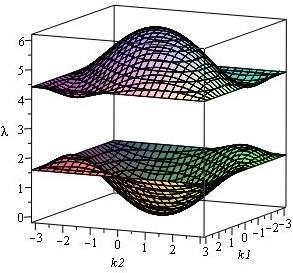}\hspace{5mm}
\includegraphics[width=4.5cm,height=4.5cm]{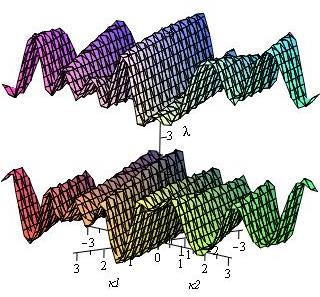}\hspace{5mm}
\includegraphics[width=4cm,height=4.5cm]{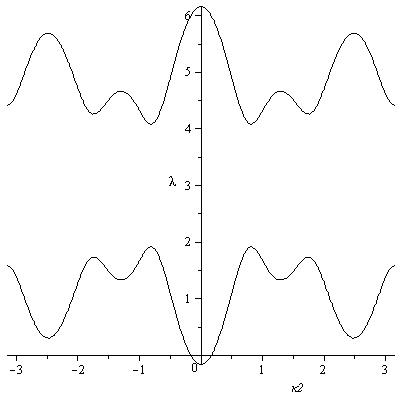}
\caption{\scriptsize Left and center: the graphs of the band functions (the dispersion relation) of the Schr\"odinger operator $H$ (when $q=1$) on the hexagonal lattice $\bG$. The Brillouin zone is parameterized by the quasimomenta $k=(k_1,k_2)$ (left) and $\k=(\k_1,\k_2)$, $\k=\wt\cT k$ (center). Right: the graphs of the band functions of the Schr\"odinger operator $H_{\gt_1}$ on the nanotube $\bG_{\gt_1}$ with the chiral vector $\gt_1=(2,3)$. The band functions of $H_{\gt_1}$ are the restrictions of the dispersion relation for $H$ (center) to the axis $\k_2$. The ranges of the corresponding band functions of $H$ and $H_{\gt_1}$ are close to each other but do not coincide.}
\label{fDR}
\end{figure}

\begin{remark}
In fact, the set $\wt\cB_\gt$ given by \er{nBzt} is a natural choice of the Brillouin zone for the subcovering $\cG_\gt$ which is a $\Z^{d-d_o}$-periodic graph with the periods $\gt_{d_o+1},\ldots,\gt_d$, see \er{neba}.
\end{remark}

\subsection{Spectra of Schr\"odinger operators on subcoverings with primitive sets of chiral vectors.}
In the proofs of our main results we will use the following known statement, see, e.g., p.421 in \cite{HJ85}.
\begin{proposition}\lb{Pmns}
For a system of linear equations $Ay=b$, where $A$ is a real $m\ts n$ matrix and $b\in\R^m$, $y\in\R^n$, the solution $\hat y$ of minimum Euclidean norm $|\hat y|$ is given by
$$
\hat y=A^+b,
$$
where $A^+$ is the Moore-Penrose pseudo-inverse of a matrix $A$. If the matrix $A$ has full row rank (i.e., $AA^\top$ is invertible), then
$$
A^+=A^\top(AA^\top)^{-1},
$$
where $A^\top$ is the transpose of $A$, and
$$
|\hat y|=|A^+b|=\big|(AA^\top)^{-1/2}b\big|.
$$
\end{proposition}

We prove Proposition \ref{PAig} about spectra of Schr\"odinger operators on subcovering graphs $\cG_\gt$ with primitive sets $\gt$ of chiral vectors.

\medskip

\no \textbf{Proof of Proposition \ref{PAig}.} Due to \er{spDt}, the spectrum of the Schr\"odinger operator $H_\gt$ on the subcovering $\cG_\gt$ is given by
\[\lb{spDt1}
\s(H_\gt)=\bigcup_{k\in\cB_\gt}\s\big(H(k)\big)=
\bigcup_{k\in\cB_\gt}\bigcup_{j=1}^\n\l_j(k)=
\bigcup_{j=1}^\n\l_j(\cB_\gt),
\]
where $\l_j(k)$ are the band functions of the Schr\"odinger operator $H$ on the periodic graph $\cG$, and the set $\cB_\gt\ss\cB:=(-\pi,\pi]^d$ is defined by \er{bgt}.  Combining \er{spDt1} and \er{wtcb}, we obtain that the spectrum of $H_\gt$ has the form \er{debe}, where the last inclusion follows from the definition \er{ban.1H} of the band edges $\l_j^\pm$ and the inclusion $\cB_\gt\ss\cB$.

\emph{i}) We prove this item for the band edge $\l_j^-(\gt)$. The proof for $\l_j^+(\gt)$ is similar. Let the condition \er{soeq} be fulfilled for some $k_o\in K_j^-$.
Then $k_o\in\cB_\gt$ (see \er{bgt}) and,
$$
\l_j^-(\gt)=\min_{k\in\cB_\gt}\l_j(k)\leq\l_j(k_o)=\l_j^-.
$$
This and the last inclusion in \er{debe} yield $\l_j^-(\gt)=\l_j^-$.

Conversely, let $\l_j^-(\gt)=\l_j^-$. Then for some $k_o\in\cB_\gt$ we have
$$
\l_j^-(\gt)=\min_{k\in\cB_\gt}\l_j(k)=\l_j(k_o)=\l_j^-.
$$
Thus, $k_o\in K_j^-$ and $k_o$ satisfies the condition \er{soeq} (since $k_o\in \cB_\gt$).

Let $0\in K_j^-$. For $k_o=0$ the condition \er{soeq} holds. Then $\l_j^-(\gt)=\l_j^-$.

Let $k_o:=(\pi,\ldots,\pi)\in K_j^-$ and the sum of the entries in each row of $\cT$ be even. Then the condition \er{soeq} holds, and $\l_j^-(\gt)=\l_j^-$.

\emph{ii}) Let $j\in\N_\n$. We prove the identity \er{as00} for $\l_j^-(\gt)$. The proof for $\l_j^+(\gt)$ is similar. Due to $2\pi\Z^d$-periodicity of the band function $\l_j(k)$, we have
$$
\l_j^-:=\min_{k\in\cB}\l_j(k)=\min_{k\in\R^d}\l_j(k)=\l_j(k_o)
$$
for some $k_o\in\wt\cB$, where the new Brillouin zone $\wt\cB$ is defined by \er{nBz}. Let $\hat k$ be the point in the set
$$
L_0=\{k\in\R^d: \cT k=0\}
$$
that is closest to the point $k_o$, i.e.,
$$
\min\limits_{k\in L_0}|k-k_o|=|\hat k-k_o|.
$$
Then $\hat k-k_o$ is the solution of the system of linear equations $\cT(k-k_o)=-\cT k_o$ with minimum Euclidean norm $|\hat k-k_o|$. Due to Proposition \ref{Pmns}, this solution has the form
$$
\hat k-k_o=-\cT^+\,\bx_o,\qqq \textrm{where}\qqq \bx_o=\cT k_o,
$$
and
\[\lb{hkmk}
|\hat k-k_o|=\big|(\cT\cT^\top)^{-1/2}\bx_o\big|.
\]
Since $k_o\in\wt\cB$, we get $\bx_o=\cT k_o\in(-\pi,\pi]^{d_o}$, see \er{nBz}. Then from \er{hkmk} we obtain
\[\lb{kokh}
\textstyle|\hat k-k_o|\leq\frac{\pi\sqrt{d_o}}\t=O\big(\frac1\t\big),
\]
where $\t$ is the smallest singular value of $\cT$. Using the last inclusion in \er{debe} and $2\pi\Z^d$-periodicity of the band function $\l_j(k)$, we have
$$
\l_j^-\leq\l_j^-(\gt)=\min_{k\in\cB_\gt}\l_j(k)=
\min_{k\in\cS_\gt}\l_j(k)\leq\l_j(\hat k),
$$
since $\hat k\in L_0\ss\cS_\gt:=\{k\in\R^d: \cT k\in2\pi\Z^{d_o}\}$. This, \er{kokh} and continuity of the band function $\l_j(\cdot)$ yield
$$
\l_j^-\leq\lim_{\t\to\iy}\l_j^-(\gt)\leq\lim_{\t\to\iy}\l_j(\hat k)
=\l_j(k_o)=\l_j^-.
$$
Thus, $\lim\limits_{\t\to\iy}\l_j^-(\gt)=\l_j^-$. \qq $\Box$

\subsection{Asymptotics of the band edges.}
We prove Theorem \ref{TAsy1} about asymptotics of the band edges $\l_j^\pm(\gt)$ of the Schr\"odinger operator $H_\gt$ on the subcovering $\cG_\gt$ as all chiral vectors from the primitive set $\gt$ are long enough.

\medskip

\no \textbf{Proof of Theorem \ref{TAsy1}.} We prove the asymptotics \er{aslk1} for $\l_j^-(\gt)$. The proof for $\l_j^+(\gt)$ is similar.  Due to \er{wtcb} and $2\pi\Z^d$-periodicity of the band function $\l_j(k)$, we have
\[\label{Tla-}
\displaystyle
\l_j^-(\gt)=\min_{k\in\wt\cB_\gt}\l_j(k)=\min_{k\in L_0}\l_j(k),
\]
where $\wt\cB_\gt$ is given by \er{nBzt} and $L_0=\big\{k\in\R^d: \cT k=0\big\}$. Thus, in order to get the band edge $\l^-_j(\gt)$ we need to find the minimum of $\l_j(k)$ with the constraint $\cT k=0$, where $\cT$ is the chiral matrix of the subcovering $\cG_\gt$ defined by \er{macT}.

Let $\tilde k_o$ be the point from the new Brillouin zone $\wt\cB$, see \er{nBz}, such that $\tilde k_o\equiv k_o\,(\hspace{-3mm}\mod 2\pi\Z^d)$. Due to $2\pi\Z^d$-periodicity, the function $\l_j(k)$ has the non-degenerate minimum $\l_j^-$ at the point $\tilde k_o$ and also satisfies the Assumption A at this point $\tilde k_o\in\wt\cB$. The Taylor expansion of the function $\l_j(k)$ about the point $\tilde k_o$ is given by
\[\lb{wtlk0}
\begin{array}{l}
\textstyle\l_j(k)=\l_j^-+\frac12\lan k-\tilde k_o,\bH(k-\tilde k_o)\ran+g_1(k-\tilde k_o),\qqq\textrm{where} \\[6pt]
\bH=\mathrm{Hess}\,\l_j(\tilde k_o)=\mathrm{Hess}\,\l_j(k_o),\\[6pt]
\l_j^-=\l_j(\tilde k_o)=\l_j(k_o),\qqq g_1(k-\tilde k_o)=O\big(|k-\tilde k_o|^3\big).
\end{array}
\]

We make the non-degenerate change of variables $k\mapsto y$, where
$$
y=\bH^{1/2}(k-\tilde k_o),
$$
and $\bH^{1/2}$ is the positive definite square root of the positive definite Hessian $\bH$. In the new variable $y$ the expansion \er{wtlk0} and the constraint $\cT k=0$ have the form
\[\lb{wtlk}
\textstyle\l_j\big(k(y)\big)=\l_j^-+\frac12|y|^2
+g_2(y),\qqq g_2(y)=g_1(\bH^{-1/2}y)=O\big(|y|^3\big),
\]
\[\lb{hupy}
\textstyle\cT(\bH^{-1/2}y+\tilde k_o)=0.
\]
Among all solutions of the system \er{hupy} we need to find the solution $\hat y$ with minimum Euclidean norm $|\hat y|$.
The $d_o\ts d$ matrix $\cT\bH^{-1/2}$ has full row rank. Thus, due to  Proposition \ref{Pmns},  the solution $\hat y$ of the system of linear equations \er{hupy} with minimum $|\hat y|$ has the form
$$
\hat y=-(\cT\bH^{-1/2})^+\,\bx_o,\qqq \bx_o=\cT\tilde k_o,
$$
and
\[\lb{nyo1}
|\hat y|^2=
\big|\big(\cT\bH^{-1/2}(\cT\bH^{-1/2})^\top\big)^{-1/2}\bx_o\big|^2=
\big|(\cT\bH^{-1}\cT^\top)^{-1/2}\bx_o\big|^2.
\]
We note that

$\bu$ $\cT\tilde k_o\in(-\pi,\pi]^{d_o}$, since $\tilde k_o\in\wt\cB$ (see \er{nBz}), and

$\bu$ $\cT\tilde k_o\equiv \cT k_o\,(\hspace{-3mm}\mod 2\pi\Z^{d_o})$, since $\tilde k_o\equiv k_o\,(\hspace{-3mm}\mod 2\pi\Z^d)$ and $\cT$ is a $d_o\ts d$ integer matrix.

\no The identity \er{nyo1} yields
\[\lb{eqno}
|\hat y|^2\leq\l_{\max}\big((\cT\bH^{-1}\cT^\top)^{-1}\big)\cdot|\bx_o|^2
=\frac{|\bx_o|^2}{\l_{\min}(\cT\bH^{-1}\cT^\top)}\,,
\]
where $\l_{\min}(A)$ and $\l_{\max}(A)$ denote, respectively, the smallest and largest eigenvalues of a self-adjoint matrix $A$. Let $\phi\in\R^{d_o}$ be the unit-norm eigenvector of $\cT\bH^{-1}\cT^\top$ corresponding to $\l_{\min}(\cT\bH^{-1}\cT^\top)$. Then
\begin{multline}\lb{lmin}
\l_{\min}(\cT\bH^{-1}\cT^\top)=\lan\phi,\cT\bH^{-1}\cT^\top\phi\ran=\lan \cT^\top\phi,\bH^{-1}\cT^\top\phi\ran\geq\l_{\min}(\bH^{-1})\lan \cT^\top\phi,\cT^\top\phi\ran\\=\l_{\min}(\bH^{-1})\lan\phi,\cT \cT^\top\phi\ran\geq\l_{\min}(\bH^{-1})\l_{\min}(\cT\cT^\top).
\end{multline}
Substituting \er{lmin} into \er{eqno} and using that $\t^2=\l_{\min}(\cT\cT^\top)$ and $\bx_o\in(-\pi,\pi]^{d_o}$, we obtain
$$
|\hat y|^2\leq\frac{\pi^2d_o}{\l_{\min}(\bH^{-1})\,\t^2}=\frac{\pi^2d_o}{\t^2}\,
\l_{\max}(\bH),
$$
which yields
\[\lb{x0gn}
|\hat y|=\textstyle O(\frac1{\t}).
\]

Using \er{Tla-}, \er{wtlk} and \er{nyo1} we get
$$
\textstyle\l_j^-(\gt)=\l_j(k(\hat y))=\l_j^-+\frac12\,|\hat y|^2
+g_2(\hat y)=\l_j^-+\frac12\,\big|(\cT\bH^{-1}\cT^\top)^{-1/2}\bx_o\big|^2
+g_2(\hat y),
$$
where $g_2(\hat y)=O(|\hat y|^3)$, or, using \er{x0gn}, $g_2(\hat y)=O(\frac1{\t^3})$. Thus, \er{aslk1} is proved.

If all components of $k_o$ are equal to either 0 or $\pi$, then $\tilde k_o\in\pi\Z^d$, and, due to $2\pi\Z^d$-periodicity and evenness of $\l_j$, we have
$$
\l_j(k+\tilde k_o)=\l_j(-k-\tilde k_o)=\l_j(-k+\tilde k_o),
$$
since $2\tilde k_o\in2\pi\Z^d$. This implies that the Taylor series \er{wtlk0} of $\l_j$ at $\tilde k_o$ has only even degree terms. Then, $g_1(k-\tilde k_o)=O\big(|k-\tilde k_o|^4\big)$ and $g_2(\hat y)=O(\frac1{\t^4})$.

Finally, let $\cG_\gt$ be a subcovering of $\cG$ with a single primitive vector $\gt\in\Z^d$. Then $\cT=\gt^\top$,
$$
\t=(\gt^\top \gt)^{1/2}=|\gt|,\qq \cT\bH^{-1}\cT^\top=\gt^\top\bH^{-1}\gt=|\bH^{-1/2}\gt|^2,\qq \bx_o=(\lan \gt,k_o\ran\hspace{-3mm}\mod2\pi)\in(-\pi,\pi],
$$
and, the asymptotics \er{aslk1} has the form \er{aslk}. \qq $\Box$

\section{Examples}
\setcounter{equation}{0}
\lb{Sec4}
In this section we consider simple examples of periodic graphs (with the unit edge weight) and describe the spectra of the Schr\"odinger operators on their subcoverings.

\subsection{Nanotubes}
The hexagonal lattice $\bG$ is a $\Z^2$-periodic graph with the periods $\ga_1,\ga_2$, see  Fig.~\ref{fig1}\emph{a}. The fundamental graph $\bG_*$ of $\bG$ consists of two vertices $v_1,v_2$ and three multiple edges connecting these vertices, see Fig.~\ref{fig1}\emph{d}.

First we briefly describe the spectrum of the Schr\"odinger operator $H=\D+Q$ on $\bG$ with a $\Z^2$-periodic potential $Q$ defined by \er{poGL}, for more details see, e.g., \cite{KP07}. The Floquet operator $H(k)$ given by \er{Flop} in the standard basis of $\ell^2(\{v_1,v_2\})\cong\C^2$ has the following matrix:
$$
H(k)=\left(
\begin{array}{cc}
  3+q & -1-e^{ik_1}-e^{ik_2} \\
  -1-e^{-ik_1}-e^{-ik_2} & 3-q
\end{array}\right), \qqq k=(k_1,k_2)\in\cB=(-\pi,\pi]^2.
$$
The eigenvalues of $H(k)$ are given by
\[\lb{evEG}
\begin{aligned}
&\l_j(k)=3+(-1)^j\sqrt{f(k)+q^2},\qq j=1,2,\qq \textrm{where}\\
&f(k)=|1+e^{ik_1}+e^{ik_2}|^2=3+2\cos k_1+2\cos k_2+2\cos(k_1-k_2).
\end{aligned}
\]
Then, using the definition of the band edges $\l_j^\pm$ in \er{ban.1H}, we obtain that the spectrum of the Schr\"odinger operator $H=\D+Q$ on the hexagonal lattice $\bG$ is given by
\[\lb{sDbG}
\begin{aligned}
&\s(H)=[\l_1^-,\l_1^+]\cup[\l_2^-,\l_2^+],\\
& \begin{array}{ll}
\l_1^-=\l_1(K_1^-)=3-\sqrt{9+q^2},\qq & \l_1^+=\l_1(K_1^+)=3-q,\\[6pt]
\l_2^-=\l_2(K_2^-)=3+q, \qq & \l_2^+=\l_2(K_2^+)=3+\sqrt{9+q^2},
\end{array}
\end{aligned}
\]
where the level sets
\[\lb{spDoK}
\textstyle K_1^-=K_2^+=\big\{(0,0)\big\}, \qqq
K_1^+=K_2^-=\big\{\big(\pm\frac{2\pi}3\,,\mp\frac{2\pi}3\big)\big\}.
\]

\medskip

Now we prove Example \ref{Ehex} about the spectrum of the Schr\"odinger operator $H_\gt=\D_\gt+Q$ on a nanotube $\bG_\gt$ with a primitive chiral vector $\gt$.

\medskip

\no \textbf{Proof of Example \ref{Ehex}.} Due to Proposition \ref{PAig}, the spectrum of the Schr\"odinger operator $H_\gt$ on the nanotube $\bG_\gt$ with a primitive chiral vector $\gt=(t_1,t_2)\in\Z^2$ consists of two bands:
$$
\s(H_\gt)=[\l_1^-(\gt),\l_1^+(\gt)]\cup[\l_2^-(\gt),\l_2^+(\gt)].
$$
Since $(0,0)\in K_1^-=K_2^+$, where $K_j^\pm$ are the level sets corresponding to the band edges $\l_j^\pm$ of the Schr\"odinger operator $H$ on the hexagonal lattice $\bG$, see \er{sDbG}, \er{spDoK}, Proposition \ref{PAig}.\emph{i} yields
$$
\l_1^-(\gt)=\l_1^-=3-\sqrt{9+q^2},\qqq \l_2^+(\gt)=\l_2^+=3+\sqrt{9+q^2}.
$$

\emph{i}) Let $t_1-t_2\in3\Z$. Then for the points $k_o=\big(\pm\frac{2\pi}3\,,\mp\frac{2\pi}3\big)\in K_1^+=K_2^-$ the condition \er{soeq} is fulfilled. Thus, due to Proposition \ref{PAig}.\emph{i} and the identities \er{sDbG},
$$
\textstyle\l_1^+(\gt)=\l_1^+=3-q,\qqq
\l_2^-(\gt)=\l_2^-=3+q,
$$
and $\ve(\gt)=0$ in \er{stin}.

\emph{ii}) Let $t_1-t_2\not\in3\Z$, i.e., $t_1-t_2=3m\pm1$ for some $m\in\Z$. We obtain the asymptotics of the band edge $\l_1^+(\gt)$ as $|\gt|\to\iy$. The proof for $\l_2^-(\gt)$ is similar. The point $k_o=\big(\frac{2\pi}3\,,-\frac{2\pi}3\big)$ is the only (up to evenness) maximum point of the band function $\l_1(k)$ in $\cB=(-\pi,\pi]^2$, see \er{sDbG}, \er{spDoK}. Using the identity \er{evEG} for $\l_1(k)$, we obtain
$$
\frac{\pa^2\l_1}{\pa k_j^2}\,(k_o)=-\frac1q\,,\qq j=1,2;
\qqq\frac{\pa^2\l_1}{\pa k_1\pa k_2}\,(k_o)=\frac1{2q}\,.
$$
Then the matrix $\bH=-\mathrm{Hess}\,\l_1(k_o)$ and its inverse are given by
$$
\bH=\frac1{2q}\left(
\begin{array}{rr}
  2 & -1 \\
  -1 & 2
\end{array}\right), \qqq \bH^{-1}=\frac{2q}3\left(
\begin{array}{cc}
  2 & 1 \\
  1 & 2
\end{array}\right).
$$
Thus, the band function $\l_1(k)$ of $H$ satisfies the Assumption A and we can use the asymptotics \er{aslk}. For $\gt=(t_1,t_2)$ and $k_o=\big(\frac{2\pi}3\,,-\frac{2\pi}3\big)$, we obtain
\[\lb{coef1}
\textstyle|\bH^{-1/2}\gt|^2=\lan \gt,\bH^{-1}\gt\ran=\frac{4q}3\,(t_1^2+t_1t_2+t_2^2),
\]
$$
\textstyle\lan\gt,k_o\ran=\frac{2\pi}3\,(t_1-t_2).
$$
Since $t_1-t_2=3m\pm1$ for some $m\in\Z$,
\[\lb{slal2}
\textstyle x_o:=(\lan\gt,k_o\ran\hspace{-3mm}\mod2\pi)=\pm\frac{2\pi}3\,.
\]
Substituting \er{coef1}, \er{slal2} and $\l_1^+=3-q$ into \er{aslk}, we obtain $$
\textstyle\l_1^+(\gt)=3-q-\frac{\pi^2}{6q(t_1^2+t_1t_2+t_2^2)}+O\big(\frac1{|\gt|^3}\big)
\qqq\textrm{as}\qqq |\gt|\to\iy.\qqq \Box
$$

\subsection{Subcoverings of the $d$-dimensional lattice}

The $d$-dimensional lattice $\dL^d$ is a $\Z^d$-periodic graph with the periods $\ga_1,\ldots,\ga_d$, for $d=3$ see  Fig.~\ref{fig2}\emph{a}. The fundamental graph $\dL^d_*$ of $\dL^d$ consists of one vertex $v$ and $d$ loop edges at this vertex $v$ (Fig.~\ref{fig2}\emph{d}). The Floquet Laplacian $\D(k)$ for the lattice $\dL^d$ has the form
\[\lb{bfDd}
\l_1(k)\equiv\D(k)=2d-2\cos k_1-\ldots-2\cos k_d, \qq k=(k_1,\ldots,k_d)\in\cB=(-\pi,\pi]^d.
\]
Thus,
$$
\textstyle \l_1^-=\min\limits_{k\in\cB}\l_1(k)=\l_1(K_1^-)=0,  \qqq \l_1^+=\max\limits_{k\in\cB}\l_1(k)=\l_1(K_1^+)=4d,
$$
where the level sets $K_1^\pm\ss\cB$ consist of the single point:
$$
K_1^-=\big\{(0,\ldots,0)\big\}, \qqq
K_1^+=\big\{\pi\1_d\big\},\qqq \1_d=(1,\ldots,1)\in\Z^d.
$$
The spectrum of the Laplacian $\D$ on the $d$-dimensional lattice $\dL^d$ is given by
$$
\s(\D)=[\l_1^-,\l_1^+]=[0,4d].
$$

Now we prove Example \ref{ExSl} about the Laplacian spectrum on a subcovering $\dL^d_\gt$ of the lattice $\dL^d$ with a primitive set $\gt$ of chiral vectors $\gt_1,\ldots,\gt_{d_o}$, where $\gt_s=(t_{s1},\ldots,t_{sd})\in\Z^d$, $s\in\N_{d_o}$, $d_o<d$.

\medskip

\no \textbf{Proof of Example \ref{ExSl}.} Due to Proposition \ref{PAig}, the spectrum of the Laplacian $\D_\gt$ on the subcovering $\dL^d_\gt$  with a primitive set $\gt$ consists of the single band:
$$
\s(\D_\gt)=[\l_1^-(\gt),\l_1^+(\gt)].
$$
Since $0\in K_1^-$, then Proposition \ref{PAig}.\emph{i} yields $\l_1^-(\gt)=\l_1^-=0$. Thus, \er{sDLt} is proved.

\emph{i}) Let $\vr_s:=t_{s1}+\ldots+t_{sd}$ be even for all $s\in\N_{d_o}$. Since $\pi\1_d\in K_1^+$, due to Proposition \ref{PAig}.\emph{i}, we have $\l_1^+(\gt)=\l_1^+=4d$.

\emph{ii}) Let $\vr_s$ be odd for some $s\in\N_{d_o}$. For the single maximum point $k_o=\pi\1_d\in K_1^+$ of the single band function $\l_1(\cdot)$ of $\D$ given by \er{bfDd} we obtain
$$
\frac{\pa^2 \l_1}{\pa k_i^2}\,(k_o)=-2,\qqq \frac{\pa^2\l_1}{\pa k_i\pa k_j}\equiv0,\qqq
i,j\in\N_d, \qqq i\neq j.
$$
The matrix $\bH=-\mathrm{Hess}\,\l_1(k_o)$ and its inverse are given by
$$
\textstyle\bH=2I_d,\qqq \bH^{-1}=\frac12\,I_d,
$$
where $I_d$ is the identity matrix of size $d$. Thus, the band function $\l_1(k)$ of $\D$ satisfies the Assumption A and we can use the asymptotics \er{aslk1}. Substituting
$$
\textstyle \l_1^+=4d,\qqq k_o=\pi\1_d, \qqq \bH^{-1}=\frac12\,I_d
$$
into \er{aslk1}, we obtain \er{ase0}.

Finally, if $\dL^d_\gt$ is a subcovering of $\dL^d$ with a single primitive chiral vector $\gt=(t_1,\ldots,t_d)\in\Z^d$ and $t_1+\ldots+t_d$ is odd, then $$
\cT=\gt^\top,\qqq \vr=(t_1+\ldots+t_d\hspace{-3mm}\mod2)=1.
$$
Substituting this into \er{ase0}, we obtain \er{ase1}. \qq $\Box$

\subsection{Subcoverings of the diamond lattice}
We consider the diamond lattice $\cG$ shown in Fig. \ref{Fcon}\emph{a}. The diamond lattice $\cG$ is obtained by stacking together infinitely many copies
$$
\bG+n\ga_3\ss\R^3,\qqq \forall\,n\in\Z,
$$
of the hexagonal lattice $\bG\ss\R^2$ along the vector $\ga_3$. Here $\R^2$ is considered as the subspace $\R^2\ts\{0\}$ of $\R^3$. The copies of $\bG$ are connected in a periodic way by the edges between white vertices in a lower copy and black vertices in the next copy of $\bG$. The diamond lattice $\cG$ is a $\Z^3$-periodic graph with the periods $\ga_1,\ga_2,\ga_3$, where $\ga_1,\ga_2$ are the periods of the hexagonal lattice $\bG$. The fundamental graph $\cG_*$ of $\cG$ consists of two vertices $v_1,v_2$ and four multiple edges connecting these vertices, see Fig.~\ref{Fcon}\emph{d}.

The spectrum of the Schr\"odinger operator $H=\D+Q$ on the diamond lattice $\cG$ with a $\Z^3$-periodic potential $Q$ defined by \er{poDL} is given by
$$
\begin{aligned}
&\s(H)=[\l_1^-,\l_1^+]\cup[\l_2^-,\l_2^+],\\
& \begin{array}{ll}
\l_1^-=\l_1(K_1^-)=4-\sqrt{16+q^2},\qq & \l_1^+=\l_1(K_1^+)=4-q,\\[6pt]
\l_2^-=\l_2(K_2^-)=4+q, \qq & \l_2^+=\l_2(K_2^+)=4+\sqrt{16+q^2},
\end{array}
\end{aligned}
$$
where the level sets $K_1^-=K_2^+=\big\{(0,0,0)\big\}$,
\[\lb{sbG3}
\textstyle  K_1^+=K_2^-=\big\{(k_1,k_2,k_3)\in\cB=(-\pi,\pi]^3: 1+e^{ik_1}+e^{ik_2}+e^{ik_3}=0\big\}.
\]
(See, e.g., Lemma 4.2 in \cite{S24}.)

\begin{figure}[t!]\centering
\includegraphics[width=6cm,height=5cm]{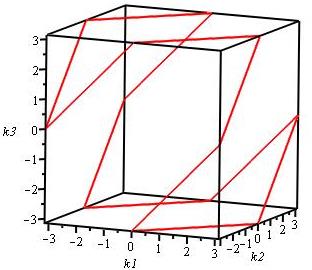}
\caption{\scriptsize The set $K_1^+=K_2^-$ defined by \er{sbG3}.}
\label{fKDl}
\end{figure}

\begin{remark}
The set $K_1^+=K_2^-$ defined by \er{sbG3} is shown in red line in Fig. \ref{fKDl}. The band edges $\l_1^+$ and $\l_2^-$ are degenerate. This example of periodic graphs where the non-degeneracy assumption fails was constructed in \cite{BK20}.
\end{remark}

Now we prove Example \ref{Exa3} describing all subcoverings $\cG_\gt$ of the diamond lattice $\cG$ isospectral to $\cG$.

\medskip

\no \textbf{Proof of Example \ref{Exa3}.} Let $\gt=(t_1,t_2,t_3)\in\Z^3$ be a primitive vector. For some $i,j=1,2,3$, $i\neq j$, the numbers $t_i$ and $t_j$ have the same parity, i.e., $t_i+t_j$ is even. Let $k_o=(k_{o1},k_{o2},k_{o3})$ be defined by
$$
k_{oi}=k_{oj}=\pi,\qqq k_{ol}=0, \qq l\in\{1,2,3\}\sm\{i,j\}.
$$
Then $k_o\in K_1^+=K_2^-$, where $K_1^+=K_2^-$ is defined by \er{sbG3} and for $k_o$ and $\cT=\gt^\top$ the condition \er{soeq} is fulfilled. For the point $0\in K_1^-=K_2^+$ the condition \er{soeq} is also fulfilled. Thus, due to Proposition \ref{PAig}.\emph{i}, the diamond lattice $\cG$ is isospectral to its subcovering $\cG_\gt$, i.e., $\s(H)=\s(H_\gt)$.

Now let $\gt=\{\gt_1,\gt_2\}\ss\Z^3$ be a primitive set, where $\gt_s=(t_{s1},t_{s2},t_{s3})$, $s=1,2$. If for some $i,j=1,2,3$, $i\neq j$, the sum $t_{si}+t_{sj}$ is even for each $s=1,2$, then using the same arguments as above we obtain that $\cG$ and $\cG_\gt$ are isospectral. Conversely, if $\cG$ and $\cG_\gt$ are isospectral, then, due to Proposition \ref{PAig}.\emph{i},
\[\lb{tk2p}
t_{s1}k_{o1}+t_{s2}k_{o2}+t_{s3}k_{o3}\in2\pi\Z,\qqq s=1,2,
\]
for some $k_o=(k_{o1},k_{o2},k_{o3})\in K_1^+=K_2^-$. Using the structure of the set $K_1^+$ (see Fig. \ref{fKDl}), we have
$$
k_{ol}=\pi,\qq k_{oi}-k_{oj}=\pm\pi,\qq \textrm{for some pairwise distinct}\;i,j,l\in\{1,2,3\}.
$$
Substituting this into \er{tk2p} and using that $t_{si}\in\Z$, we obtain
\[\lb{tk3p}
(t_{si}+t_{sj})k_{oj}+\pi(t_{si}+t_{sl})\in2\pi\Z,\qqq s=1,2.
\]
We assume that there are no $i,j=1,2,3$, $i\neq j$, such that the sum $t_{si}+t_{sj}$ is even for each $s=1,2$, i.e., for each pair $(i,j)\in\{1,2,3\}^2$ with $i\neq j$, at least one of the sums $t_{si}+t_{sj}$, $s=1,2$, is odd. This yields that, up to the permutations of the index $s=1,2$, only the following three cases are possible:
\begin{table}[h]
\begin{tabular}{|c|c|c|c|c|}
  \hline
Case & $t_{1i}+t_{1j}$ & $t_{2i}+t_{2j}$  & $t_{1i}+t_{1l}$ & $t_{2i}+t_{2l}$\\ \hline
 1 & even & odd  & odd & even \\ \hline
 2 & even & odd  & odd & odd \\ \hline
 3 & odd & odd  & even & odd \\ \hline
\end{tabular}
\end{table}

We consider Case 1. Cases 2 and 3 are similar. In Case 1 the conditions \er{tk3p} have the form
$$
2n_1k_{oj}=(2m_1+1)\pi,\qqq (2n_2+1)k_{oj}=2m_2\pi
$$
for some $n_s,m_s\in\Z$, $s=1,2$, which yields the contradiction $4n_1m_2=(2m_1+1) (2n_2+1)$. Thus, if $\cG$ and $\cG_\gt$ are isospectral, then for some $i,j=1,2,3$, $i\neq j$, the sum $t_{si}+t_{sj}$ is even for each $s=1,2$.  \qq $\Box$

\medskip

\textbf{Acknowledgments.}  This work is supported by the Russian Science Foundation (project No. 25-21-00157). I would like to thank the referees for thoughtful comments that helped me to improve the manuscript.

\medskip

\textbf{Data availability}  No datasets were generated or analysed during the current study.

\medskip

\textbf{Conflict of interest} None declared.



\medskip

\begin{thebibliography}{9999}
\setlength{\itemsep}{-\parskip}\footnotesize

\bibitem[A95]{A95} Adachi, T. On the spectrum of periodic Schr\"odinger operators and a tower of coverings. Bull. London Math. Soc. 27 (1995), no. 2, 173--176.

\bibitem[AM76]{AM76} Ashcroft, N.W.; Mermin, N.D. Solid State Physics, Holt, Rinehart and Winston, New York -- London, 1976.

\bibitem[BCCM22]{BCCM22} Berkolaiko, G.; Canzani, Y.; Cox, G.; Marzuola, J.L. A local test for global extrema in the dispersion relation of a periodic graph. Pure Appl. Anal. 4 (2022), no. 2, 257--286.

\bibitem[BK13]{BK13} Berkolaiko, G.; Kuchment, P. Introduction  to
Quantum Graphs, Mathematical Surveys and Monographs, V. 186 AMS, 2013.

\bibitem[BK20]{BK20} Berkolaiko, G.; Kha, M. Degenerate band edges in periodic quantum graphs. Lett. Math. Phys. 110 (2020), no. 11, 2965--2982.

\bibitem[C97]{C97}  Cassels, J.W.S. An introduction to the geometry of numbers, Springer, Berlin, 1997.

\bibitem[D15]{D15} Do, N. On the quantum graph spectra of graphyne nanotubes. Anal. Math. Phys. 5 (2015), 39--65.

\bibitem[GKT93]{GKT93} Gieseker, D.; Kn\"orrer, H.; Trubowitz, E. The geometry of algebraic Fermi curves, volume 14 of Perspectives in Mathematics. Academic Press Inc., Boston, MA, 1993.

\bibitem[FS24]{FS24} Faust, M.; Sottile, F. Critical points of discrete periodic operators, J. Spectr. Theory 14, (2024), no. 2,  1--35.

\bibitem[HKSW07]{HKSW07} Harrison, J.M.; Kuchment, P.; Sobolev, A.; Winn, B.  On occurrence of spectral edges for periodic operators inside the Brillouin zone, J. Phys. A 40 (2007), 7597--7618.

\bibitem[HN09]{HN09} Higuchi, Y.; Nomura, Y. Spectral structure of the Laplacian on a covering graph, Eur. J. Combin. 30 (2009), no. 2, 570--585.

\bibitem[HS99]{HS99} Higuchi, Y.; Shirai, T. A remark on the spectrum of magnetic Laplacian on a graph, the proceedings of TGT10, Yokohama Math. J., 47 (Special issue) (1999), 129--142.

\bibitem[HS04]{HS04} Higuchi, Y.; Shirai, T. Some spectral and geometric properties for infinite graphs, AMS Contemp. Math. 347 (2004), 29--56.

\bibitem[HJ85]{HJ85} Horn, R.; Johnson, C. Matrix analysis. Cambridge University Press, 1985.

\bibitem[K89]{K89} Kappeler, T. Isospectral potentials on a discrete lattice. III, Trans. Amer. Math. Soc. 314 (1989), no. 2, 815--824.

\bibitem[KKR17]{KKR17} Kha, M.; Kuchment, P.; Raich, A. Green's function asymptotics near the internal edges of spectra of periodic elliptic operators. Spectral gap interior. J. Spectr. Theory 7 (2017), no.4, 1171--1233.

\bibitem[KS21]{KS21} Koll\'{a}r, A.J.; Sarnak, P. Gap sets for the spectra of cubic graphs, Comm. Amer. Math. Soc. 1 (2021), 1--38.

\bibitem[KK10]{KK10} Korotyaev, E.;  Kutsenko, A. Zigzag and armchair nanotubes in external fields, Adv. Math. Res. 10 (2010), 273--302.

\bibitem[K05]{K05}  Kuchment, P. Quantum graphs. II. Some spectral properties of quantum and combinatorial graphs. J. Phys. A 38 (2005), no. 22, 4887--4900.

\bibitem[KP07]{KP07} Kuchment, P.; Post., O. On the spectra of carbon nano-structures, Comm. Math. Phys. 275 (2007), no. 3, 805--826.

\bibitem[L23]{L23} Liu, W. Floquet isospectrality for periodic graph operators, J. Differ. Equ. 374 (2023), 642--653.

\bibitem[L24]{L24} Liu, W. Fermi isospectrality for discrete periodic Schr\"odinger operators, Comm. Pure Appl. Math. 77 (2024), 1126--1146.

\bibitem[RS78]{RS78} Reed, M.; Simon, B. Methods of modern mathematical
physics, vol.IV. Analysis of operators, Academic Press, New York, 1978.

\bibitem[SY23]{SY23} Sabri, M.; Youssef, P. Flat bands of periodic graphs, J. Math. Phys. 64 (2023), no. 9.

\bibitem[S24]{S24} Saburova, N. Asymptotic isospectrality of Schr\"odinger operators on periodic graphs, Anal. Math. Phys. 14 (2024), no. 74.

\bibitem[S13]{S13} Sunada, T. Topological crystallography, Surveys Tutorials  Appl. Math. Sci., vol. 6, Springer, Tokyo, 2013.

\bibitem[SS92]{SS92} Sy, P.W.; Sunada, T. Discrete Schr\"odinger
 operator on a graph, Nagoya Math. J. 125 (1992), 141--150.

\end{thebibliography}
\end{document}